\documentclass[reqno, 12pt]{amsart}
\usepackage[all]{xy}
 \RequirePackage{amsgen}
\RequirePackage{amstext} \RequirePackage{amsbsy}
\RequirePackage{amsopn} \RequirePackage{amsthm}
\RequirePackage{amssymb} \RequirePackage{epsfig}
\RequirePackage{amscd}

\usepackage[mathscr]{eucal}

\newtheorem{lemma}{Lemma}[section]
\newtheorem{cor}{Corollary}[section]
\newtheorem{prop}{Proposition}[section]

\newtheorem{theorem}{Theorem}[section]

\theoremstyle{definition}
\newtheorem{defn}{Definition}[section]
\newtheorem{conj}{Conjecture}[section]

\theoremstyle{remark}
\newtheorem{rem}{Remark}[section]

\let\a\alpha  \let\b\beta  \let\c\gamma  \let\d\delta
  
  \let\l\lambda   
      
\let\GL\Lambda
\def\C{\mathbb C}

\def\CH{\mathcal H}
\def\CO{\mathcal O}
\def\CL{\mathcal L}
\def\E{\mathcal E}
\def\CF{\mathcal F}

\def\CL{\mathcal L}
\def\CS{\mathcal S}

\def\tG{\widetilde G}

\def\sgn{\rm sgn}

\def\mG{\rm{G}}
\def\GL{\rm{GL}}
\def\p{{\mathfrak p}}
\def\Tr{{\rm Tr}}
\def\Gal{{\rm Gal}}
\def\Frob{{\rm Frob}}
\def\F{{\mathbb F}}
  
\def\Z{{\mathbb Z}}

\def\Q{{\mathbb Q}}

\def\O{\mathcal O}

\def\P{\wp}

\def\la{\langle}
\def\ra{\rangle}
\def\G{\Gamma}
\def\ord{{\rm ord}}

\title{ON ATKIN-SWINNERTON-DYER CONGRUENCE RELATIONS}

\author{Wen-Ching Winnie Li}

\address{Department of Mathematics \\Pennsylvania State University\\
University Park, PA 16802\\USA} \email{wli@math.psu.edu}

\author{Ling Long}
\address{Department of Mathematics\\Iowa State University\\ Ames, IA 50011
\\USA}
\email{linglong@iastate.edu}

\author{Zifeng Yang}
\address{Department of Mathematics\\Capital Normal University\\
 Beijing 100037\\P. R. China}
\email{yangzf@mail.cnu.edu.cn}
\date{}

\thanks{The research of the first author was supported in part by a NSF grant
 DMS 99-70651 and a NSA grant MDA904-03-1-0069.
The second author was supported in part by a NSF grant No. DMS
97-29992 and a Liftoff grant from the Clay Mathematical
Institute.}

\begin{document}

\begin{abstract}
In this paper we exhibit a noncongruence subgroup $\G$ whose space
of weight 3 cusp forms $S_3(\G)$ admits a basis satisfying the
Atkin-Swinnerton-Dyer congruence relations with two weight 3
newforms for certain congruence subgroups. This gives a modularity
interpretation of the motive attached to $S_3(\G)$ by A. Scholl
and also verifies the Atkin-Swinnerton-Dyer congruence conjecture
for this space.

\end{abstract}

\maketitle
\section{Introduction}\label{Introduction}

The theory of modular forms for congruence subgroups is well
developed. Given a cuspidal normalized newform $g = \sum_{n \ge 1}
a_n(g)q^n$, where $q = e^{2\pi i\tau}$, of weight $k \ge 2$ level
$N$ and character $\chi$, the Fourier coefficients of $g$ satisfy
the recursive relation
\begin{eqnarray}\label{e:0.1}
 a_{np}(g) - a_p(g)a_n(g) + \chi(p)p^{k-1}a_{n/p}(g) = 0
\end{eqnarray}
for all primes $p$ not dividing $N$ and for all $n \ge 1$. Thanks
to the work of Eichler, Shimura, and Deligne
\cite{eich57,shim59,del1}, there exists a compatible family of
$\lambda$-adic representations $\rho_{\lambda,g}$ of the Galois
group $\Gal(\bar\Q/\Q)$, unramified outside $lN$, where $\lambda$
divides $l$, such that
\begin{eqnarray*}\Tr(\rho_{\lambda,g}(\Frob_p)) &=& a_p(g), \\
\det(\rho_{\lambda,g}(\Frob_p)) &=& \chi(p)p^{k-1},
\end{eqnarray*}
for all primes $p$ not dividing $lN$. Combining both, we see that
the characteristic polynomial $H_p(T) = T^2 - A_1(p)T + A_2(p)$ of
$\rho_{\lambda, g}(\Frob_p)$ is independent of the $\lambda$'s not
dividing $p$, and the Fourier coefficients of $g$ satisfy the
relation
\begin{eqnarray}\label{e:0.2}
 a_{np}(g) - A_1(p)a_n(g) + A_2(p)a_{n/p}(g) = 0
\end{eqnarray}
for all $n \ge 1$ and all primes $p$ not dividing $N$.

The knowledge on modular forms for noncongruence subgroups,
however, is far from satisfactory. For example, the paper
\cite{Sch97} by Scholl and Serre's letter to Thompson in
\cite{Th89} indicate that the obvious definitions
 of Hecke operators for noncongruence subgroups would not work well. On the
positive side, Atkin and Swinnerton-Dyer \cite{ASD} initiated the
study of the arithmetic properties of modular forms for
noncongruence subgroups, and they have observed very interesting
congruence relations for such forms, which we now explain.

Let $\G$ be a noncongruence subgroup of $SL_2(\Z)$ with finite
index. For an integer $k \ge 2$, denote by $S_k(\G)$ the space of
cusp forms of weight $k$ for $\G$ and by $d$ its dimension. The
modular curve $X_{\G}$, the compactification by adding cusps of
the quotient of the Poincar\'e upper half plane by $\G$, has a
model defined over a number field $K$ in the sense of Scholl
\cite[\S 5]{Sch85}.

 As explained in \cite{ASD,Sch85,Sch87},
there exists a subfield $L$ of $K$, an element $\kappa \in K$ with
$\kappa^\mu \in L$, where $\mu$ is the width of the cusp
$\infty$, and a positive integer $M$ such that $\kappa^\mu$ is
integral outside $M$ and $S_k(\G)$ has a basis consisting of
$M$-integral forms. Here a form $f$ of $\G$ is called $M$-integral
if in its Fourier expansion at the cusp $\infty$

\begin{eqnarray}\label{e:0.4}
f(\tau) = \sum_{n \ge 1} a_n(f)q^{n/\mu},
\end{eqnarray}
the Fourier coefficients $a_n(f)$ can be written as $\kappa^n
c_n(f)$ with $c_n(f)$ lying in the ring $\CO_L[1/M]$, where
$\CO_L$ denotes the ring of integers of $L$.
\medskip

Based on their numerical data, Atkin and Swinnerton-Dyer
\cite{ASD} made an amazing discovery of congruence relations for
certain cusp forms for noncongruence subgroups. It is tempting to
extrapolate from their observations and from our own numerical
data to formulate the following congruence conjecture.

\begin{conj} [Atkin-Swinnerton-Dyer congruences] Suppose that the modular curve
$X_{\G}$ has a model over $\Q$ in the sense of \cite[\S 5]{Sch85}.
There exist a positive integer $M$ and a basis of $S_k(\G)$
 consisting of $M$-integral forms $f_j$, $1 \le j \le d$, such that for each
prime $p$ not dividing $M$, there exists a nonsingular $d \times
d$ matrix $(\lambda_{i,j})$ whose entries are in
a finite extension of $\Q_p$, algebraic integers $A_p(j)$, $1 \le
j \le d$, with $|\sigma(A_p(j))| \le 2 p^{(k-1)/2}$ for all
embeddings $\sigma$, and characters $\chi_j$ unramified outside
$M$ so that for each $j$ the Fourier coefficients of $h_j : =
\sum_{i}\lambda _{i,j} f_i$  satisfy the congruence relation
\begin{eqnarray}\label{e:0.5}
\ord_p (a_{np}(h_j) - A_p(j)a_n(h_j) +
\chi_j(p)p^{k-1}a_{n/p}(h_j)) \ge (k - 1)(1 + \ord_pn)
\end{eqnarray}
for all $n \ge 1$; or equivalently, for all $n \ge 1$,
\begin{eqnarray*}
(a_{np}(h_j) - A_p(j)a_n(h_j) +
\chi_j(p)p^{k-1}a_{n/p}(h_j))/(np)^{k-1}
\end{eqnarray*}
is integral at all places dividing $p$.
\end{conj}

In other words, the recursive relation (1) on Fourier coefficients
of modular forms for congruence subgroups is replaced by the
congruence relation (\ref{e:0.5}) for forms of noncongruence
subgroups. The meaning of $A_p(j)$'s is mysterious; the examples
in \cite{ASD} suggest that they satisfy the Sato-Tate conjecture.
\medskip

In \cite{Sch85} Scholl proved a ``collective version" of this
conjecture.
\medskip

\begin{theorem}[Scholl] Suppose that $X_{\G}$ has a model over $\Q$ as before.
Attached to $S_k(\G)$ is a compatible family of $2d$-dimensional
$l$-adic representations $\rho_l$ of the Galois group $\Gal(\bar
\Q/\Q)$ unramified outside $lM$
 such that for primes $p > k + 1$ not dividing $Ml$, the following hold.

(i) The characteristic polynomial
\begin{eqnarray}\label{e:0.6}
H_p(T) = \sum_{0 \le r \le 2d} B_r(p)T^{2d-r}
\end{eqnarray}
of $\rho_l(\Frob_p)$ lies in $\Z[T]$ and is independent of $l$,
and its roots
 are algebraic integers with absolute value $p^{(k-1)/2}$;

(ii) For any $M$-integral form $f$ in $S_k(\G)$, its Fourier
coefficients $a_n(f)$, $n \ge 1$, satisfy the congruence relation
\begin{equation}\label{e:0.7}
\begin{split}
&\ord_p(a_{np^d}(f) + B_1(p)a_{np^{d-1}}(f) + \cdots +
B_{2d-1}(p)a_{n/p^{d-1}}(f) + B_{2d}(p)a_{n/p^d}(f)) \\
&\ge (k-1)(1 + \ord_pn)
\end{split}
\end{equation}
for $n \ge 1$.
\end{theorem}

Scholl's theorem establishes the Atkin-Swinnerton-Dyer congruences
 if $S_k(\G)$ has dimension 1. If the Atkin-Swinnerton-Dyer congruences
were established in general, then
\begin{eqnarray*}
H_p(T) = \prod_{1 \le j \le d}(T^2 - A_p(j)T + \chi_j(p)p^{k-1}).
\end{eqnarray*}
Scholl's congruence relation (\ref{e:0.7}) may be regarded as a
collective replacement for forms of noncongruence subgroup of the
equality (\ref{e:0.2}) for newforms.

\medskip

Let $f=\sum_{n \ge 1}a_n(f)q^{n/\mu}$ be a $M$-integral cusp form
in $S_k(\G)$, and let $g = \sum_{n \ge 1}b_n(g)q^n$ be a
normalized newform of weight $k$ level $N$ and character $\chi$.

\begin{defn} The two forms $f$ and $g$ above are said to satisfy the
Atkin-Swinnerton-Dyer congruence relation if, for all primes $p$
not dividing $MN$ and for all $n \ge 1$,
\begin{eqnarray}\label{e:0.8}
 (a_{np}(f) - b_p(g)a_n(f) + \chi(p)p^{k-1}a_{n/p}(f))/(np)^{k-1}
 \end{eqnarray}
is integral at all places dividing $p$.
\end{defn}

In particular, if $S_3(\G)$ has a basis of $M$-integral forms
$f_j, 1 \le j \le d$, such that each $f_j$ satisfies the
Atkin-Swinnerton-Dyer congruence relation with some cuspidal
newform $g_j$ of weight 3 for certain congruence subgroup, then
this not only establishes the Atkin-Swinnerton-Dyer congruences
conjecture for the space $S_3(\G)$, but also provides an
interpretation of the $A_p(j)$'s in the conjecture. Geometrically,
this means that the motive attached to $S_3(\G)$ by Scholl comes
from modular forms for congruence subgroups. Furthermore, if $-I
\notin \G$ and $\E_{\G}$ is the elliptic modular surface
associated to $\G$ in the sense of \cite{Shi72}, then $\E_{\G}$ is
an elliptic surface with base curve  $X_{\G}$. The product of the
$L$-functions $\prod_{1 \le j \le d}L(s, g_j)$ occurs in the
Hasse-Weil $L$-function $L(s, \E_{\G})$ attached to the surface
$\E_{\G}$, and it is the part arising from the transcendental
lattice of the surface. In this case, $L(s, \E_{\G})$ has both its
numerator and denominator product of automorphic $L$-functions. In
other words, the Hasse-Weil $L$-function of $\E_{\G}$ is
``modular". To-date, only a few such examples are known. In a
recent preprint of Livn\'e and Yui \cite{l-y}, the $L$-functions
of some rank 4 motives associated to non-rigid Calabi-Yau
threefolds are proven to be the $L$-functions of some automorphic
forms.
\smallskip

For the noncongruence subgroup $\G_{7,1,1}$ studied in \cite{ASD},
the space $S_4(\G_{7,1,1})$ is one-dimensional. Let $f$ be a
nonzero 14-integral form in $S_4(\G_{7,1,1})$. Scholl proved in
\cite{Sch88} that there is a normalized newform $g$ of weight 4
level 14 and trivial character such that $f$ and $g$ satisfy the
Atkin-Swinnerton-Dyer congruence relation. In the unpublished
paper \cite{Sch93}, Scholl obtained a similar result for
$S_4(\G_{4,3})$ and $S_4(\G_{5,2})$; both spaces are also
1-dimensional.
\smallskip

The purpose of this paper is to present an example of
2-dimensional space of cusp forms of weight 3 whose associated $l$-adic
representation is modular and the existence of
a $M$-integral basis, independent of $p$, such that each satisfies
the Atkin-Swinnerton-Dyer congruence relation with a cusp form of
a congruence subgroup. More precisely, we shall prove

\medskip

\begin{theorem} Let $\G$ be the index 3 noncongruence subgroup of $\G^1(5)$ such that
the widths at two cusps $\infty$ and $-2$ are 15.
\begin{itemize}
\item[(1)] Then $X_{\G}$
has a model over $\Q$, $\kappa = 1$, and the space $S_3(\G)$ is
2-dimensional with a basis  consisting of $3$-integral forms
\begin{eqnarray*}
f_+(\tau) &=&  q^{1/15} + iq^{2/15} - \frac {11}{3}q^{4/15}
-i\frac{16}{3}q^{5/15} -\frac {4}{9}q^{7/15} + i\frac
{71}{9}q^{8/15} + \frac{932}{81}q^{10/15}
+ O(q^{11/15}),  \\
f_-(\tau) &=& q^{1/15} - iq^{2/15} - \frac {11}{3}q^{4/15}
+i\frac{16}{3}q^{5/15} -\frac {4}{9}q^{7/15} - i\frac
{71}{9}q^{8/15} + \frac{932}{81}q^{10/15} + O(q^{11/15}).
\end{eqnarray*}
\item[(2)] The $4$-dimensional $l$-adic representation $\rho_{l}$ of 
$\Gal (\bar \Q/\Q)$ associated to $S_3(\G)$ constructed by
Scholl is modular. More precisely, there are two cuspidal newforms of weight 3 level 27 and
character $\chi_{-3}$ given by
\begin{eqnarray*}
g_+(\tau) & = & q - 3iq^2 - 5q^4 + 3iq^5 + 5q^7 + 3iq^8 + 9q^{10}
+ 15iq^{11} - 10q^{13} - 15iq^{14} - 11q^{16} \\&& - 18iq^{17} -
16q^{19} - 15iq^{20} + 45q^{22} + 12iq^{23} + O(q^{24})\\
g_-(\tau) & = & q + 3iq^2 - 5q^4 - 3iq^5 + 5q^7 - 3iq^8 + 9q^{10}
- 15iq^{11} - 10q^{13} + 15iq^{14} - 11q^{16} \\&& + 18iq^{17} -
16q^{19} + 15iq^{20} +
45q^{22} - 12iq^{23} + O(q^{24})\\
\end{eqnarray*}
such that over the extension by joining $\sqrt {-1}$, $\rho_{l}$
decomposes into the direct sum of the two $l$-adic representations attached to
$g_+$ and $g_{-}$.

\item[(3)] $f_{+}$ and $g_{+}$ (resp. $f_{-}$ and $g_{-}$) satisfy
the Atkin-Swinnerton-Dyer congruence relation.
\end{itemize}
\end{theorem}

\noindent Here $\chi_{-3}$ is the quadratic character attached to
the field $\Q(\sqrt {-3})$. The precise definition of $\G$ in
terms of generators and relations is given at the end of \S3.
\medskip

The proof of this theorem occupies \S2 - \S7. Here we give a
sketch. The modular curve $X_{\G}$ of $\G$ is a three fold cover
of the congruence modular curve $X_{\G^1(5)}$ ramified only at two
cusps of $\G^1(5)$. By explicitly computing the Eisenstein series
of weight 3 for $\G^1(5)$, we obtain in \S4 an explicit basis
$f_+$ and $f_-$ of $S_3(\G)$ which are $3$-integral, as stated
above.

To establish the congruence relations, we take advantage of the
existence of an elliptic surface $\E$ over $X_{\G}$ with an
explicit algebraic model. There exists a $\Q$-rational involution
$A$ on $X_{\G}$ which induces an action on $\E$ of order 4, which
commutes with the action of the Galois group over $\Q$. In fact,
$f_+$ and $f_-$ are eigenfunctions of $A$ with eigenvalues $-i$
and $i$, respectively. The explicit defining equation of $\E$
gives rise to a 4-dimensional $l$-adic representation $\rho_l^*$
of $\Gal(\bar \Q/\Q)$ which is isomorphic to the $l$-adic
representation $\rho_l$ Scholl constructed in \cite{Sch85} at most
up to a quadratic twist $\phi$. Take $l= 2$. The dyadic
representations are unramified outside 2 and 3.  Making use of the
action of $A$, we may regard $\rho_2^*$ as a two-dimensional
representation over $\Q_2(i)$, which is isomorphic to the
completion of $\Q(i)$ at the place with $1+i$ as a uniformizer,
denoted by $\Q(i)_{1+i}$ for convenience. The explicit defining
equation allows us to determine the characteristic polynomial of
the $\Frob_p$ under $\rho_2^*$ over $\Q_2$ for small primes, and
that over $\Q(i)_{1+i}$ except for primes congruent to 2 mod 3, in
which case the trace is determined up to sign (cf. Table (1)).

On the other hand, the two cuspidal newforms $g_+$ and $g_-$
combined come from a 4-dimensional 2-adic representation $\tilde
\rho_2$ of $\Gal(\bar \Q/\Q)$, on whose space the Atkin-Lehner
operator $H_{27}$ acts. It has order 4. So we may also regard
$\tilde \rho_2$ as a 2-dimensional representation over
$\Q(i)_{1+i}$. The characteristic polynomials of the Frobenius
elements under $\tilde \rho_2$ are easily read off from the
Fourier coefficients of $g_{\pm}$, which we obtained from W.
Stein's website \cite{stein}. Since the residue field of
$\Q(i)_{1+i}$ is $\F_2$, we use Serre's method \cite{Serre} to
show that $\rho_2^*$ and $\tilde \rho_2$ are isomorphic from the
incomplete information of the characteristic polynomials of the
Frobenius elements at primes $5 \le p \le 19$.
This also implies that $\rho_2^*$ is isomorphic to its twist by
the quadratic character $\chi_{-3}$.

To prove that $\rho_2^*$ and $\rho_2$ are isomorphic, we show that
$\rho_2^*$ satisfies the congruence relations (\ref{e:0.7}) for $n
= 1$, $p = 7$ and $p = 13$. This in turn forces $\phi$ to be
either $\chi_{-3}$ or trivial. In either case, we have the desired
isomorphism. Finally, to obtain the Atkin-Swinnerton-Dyer
congruence relations between $f_{\pm}$ and $g_{\pm}$,
 we compare the $p$-adic theory on eigenspaces of $A$ and the dyadic theory on
eigenspaces of $H_{27}$, and draw the desired conclusion using Scholl's proof
of Theorem 1.1 in \S5 of \cite{Sch85}.

We end the paper by observing that if the space of cusp forms of
weight 3 for a noncongruence subgroup $\G'$ is 1-dimensional with
a nonzero $M$-integral form $f$ and there is an elliptic $K3$
surface over the modular curve $X_{\G'}$, then there is a cuspidal
newform $g$ such that $f$ and $g$ satisfy the
Atkin-Swinnerton-Dyer congruence relations.

The authors would like to thank Prof. J.-P. Serre for many
stimulating and helpful communications. We are particularly
grateful to him for explaining how to apply his method to our
situation. Special thanks are also due to Prof. W. Hoffman for his
numerous suggestions and thought-provoking questions which led to
substantial improvements of the paper.

\bigskip

\section{An elliptic surface}\label{some_elliptic_surface}
Let $\E$ denote the minimal smooth model of the elliptic surface
given by
\begin{equation}\label{e:1.1}
  y^2+(1-t^3)xy-t^3y=x^3-t^3x^2\,,
\end{equation} where the parameter $t$ runs through the points in the
complex projective line $\C P^1$. Viewed as an elliptic curve
defined over $\C(t)$, its $j$-invariant is
\begin{equation}\label{j}
j=\frac{(t^{12}-12t^9+14t^6+12t^3+1)^3}{t^{15}(t^6-11t^3-1)}\,.
\end{equation} Its Mordell-Weil group is isomorphic to $\Z/5\Z$.
Indeed, it is a subgroup of $\Z/15\Z \times \Z/15\Z$ and it
contains $\Z/5\Z$ as a subgroup. By examining the restriction of
the determinant of its transcendental lattice, we find that this
group has order a power of 5. We conclude that the group is
$\Z/5\Z$ from a result of D. A. Cox and W. R. Parry \cite{c-p}.

It is clear that at a generic $t\in \C$, the fiber of the natural
projection
\begin{eqnarray}\label{f}
\pi:\ E & \rightarrow & \C P^1 \\
(x,y,t)& \mapsto & t
\end{eqnarray}
\noindent is an elliptic curve, namely, a smooth compact curve of
genus one. This is the case except for 8 values of $t$: 0,
$\infty$, and the six roots of $t^6-11t^3-1=0$. At these 8
exceptional values of $t$, there are 8 special fibers.  They are
identified, using Tate algorithm, to be of respective type
$I_{15}$, $I_{15}$, $I_1$, $I_1$, $I_1$, $I_1$, $I_1$, $I_1$ in
Kodaira's notation. The surface $\E$ is an elliptic modular
surface by [Shi72], [Nor85]. Denote by $\Gamma$ an associated
modular group, which is a subgroup of $SL_2(\Z)$ of finite index.
Let ${\bar \Gamma}=\pm \Gamma/\pm I$ be its projection in
$PSL_2(\Z)$.  We know from the information of the special fibers
that the group $\Gamma$ has no elliptic points. In another words,
it is a torsion free subgroup of $SL_2(\Z)$.

\smallskip

Now we consider some topological and geometrical invariants of
$\E$. By Kodaira's formula [Kod63], the Euler characteristic of
$\E$ is 36. Its irregularity, which equals the genus of the base
curve, is 0. Its geometric genus is 2 by Noether's formula. Denote
by  $h^{i,j}=\dim H^j(\E, \Omega_{\E}^i)$ the $(i,j)$'s Hodge
number of $\E$. Then the Hodge numbers of $\E$ can be arranged
into the following Hodge diamond
\[
\begin{split}
&1\\
0\quad &\ \quad 0\\
2\quad\quad &30\quad\quad 2\\
0\quad &\ \quad 0\\
&1
\end{split}
\] where on the $(i+1)$th row, the $(j+1)$th number is
$h^{i-j,j},$ with $j=0,..., i$ for $0 \le i \le 2$,  and $j=0,4-i$
for $3 \le i \le 4$.

As the group $\bar \Gamma$ is not among the genus-zero torsion
free congruence subgroups of $PSL_2(\Z)$ listed by A. Sebbar
[Seb01], we conclude that $\Gamma$ is a noncongruence subgroup of
$SL_2(\Z)$.

 Let $L$ be the free part of the cohomology
group $H^2(\E,\Z)$.  It is an even unimodular lattice with the
bilinear form given by the cup-product. The signature of this
lattice is $(5,29)$ by the Hodge index theory. It follows from the
classification of even unimodular lattices that $L$ is isometric
to $U^5\oplus E_8(-1)^3$, where $U$ denotes the hyperbolic matrix
and $E_8(-1)$ denotes the unique negative definite even unimodular
lattice of rank 8.

By the Shioda-Tate formula [Shi72], the Picard number is $
2+2(15-1)=30. $ The N\'eron-Severi group $NS(\E)$, which is the
group of divisors on $\E$ modulo algebraic equivalence, is a
torsion-free $\Z$-module of rank 30. This group can be imbedded
into $L$ by a cohomology sequence. The determinant of this
sublattice, by a formula in [Shi72], is equal to
\[
|\det(NS(\E))|=\frac{15^2}{5^2}=9.
\]
The orthogonal complement $T_{\E}$ of $NS(\E)$ in $L$, called the
transcendental lattice of $L$, has rank 4 and
$|\det(T_{\E})|=|\det(NS(\E))| = 9$  since $L$ is unimodular,

As we are interested in the arithmetic properties of $\E$, we
shall consider the reductions of $\E$. It turns out that for this
particular elliptic surface $\E$ the only bad prime is 3. (The
prime 5 is good because the 5 torsion points have killed the
contribution of 5 from the special fibres.) Hence we may regard
$\E$ as a normal connected smooth scheme  over $\Z[1/3]$; it is
tamely ramified along the closed subscheme formed by the cusps.

\section{Determining the noncongruence subgroup
$\Gamma$}\label{construction} For any positive integer $N$ let
\[
\begin{split}
&\Gamma^0(N)=\left\{ \left( \begin{matrix}
             a & b\\
             c & d
             \end{matrix}
            \right)
        \in SL_2(\Z): N|b \right\}\,, \\
&\Gamma^1(N)=\left\{ \left( \begin{matrix}
             a & b\\
             c & d
              \end{matrix}
            \right)
        \in \Gamma^0(N): a\equiv  d \equiv 1 \mod N \right\}\,.
\end{split}
\]
With $t^3$ in (\ref{e:1.1}) replaced by $t$, the new equation
defines an elliptic modular surface $\E'$ over the modular curve
for the group $\Gamma^1(5)$. The surface $\E$ is a three fold
cover of $\E'$, and thus the group $\Gamma$ is a subgroup of
$\Gamma^1(5)$ of index three. We proceed to determine $\Gamma$ in
terms of generators and relations.

First we decompose the full modular group $SL_2(\Z)$ as
\[
SL_2(\Z) = \bigcup_{1 \le i \le 6} \Gamma^0(5) \gamma_i,
\]
where $\gamma_i = \displaystyle{\left(\begin{matrix}1 &i-1\\ 0 &
1\end{matrix}\right)}$ for $1 \le i \le 5$ and $\gamma_6 =
\displaystyle{\left(\begin{matrix}0 & -1\\1 &
0\end{matrix}\right)}$. Further,
$$
\Gamma^0(5)=\bigcup_{1 \le j \le 4} \Gamma^1(5)A^j =
\pm\Gamma^1(5)\bigcup\pm\Gamma^1(5)A,
$$
\noindent where $A=\displaystyle{\left(\begin{matrix}-2 & -5\\1 &
2\end{matrix} \right)}$, $A^2=-I$. Hence the coset representatives
of $\pm\Gamma^1(5)$ in $SL_2(\Z)$ may be taken as $\gamma_i$ and
$A\gamma_i$ for $1\le i \le 6$.

Listed below are the cusps of $\pm\Gamma^1(5)$ and a choice of
generators of their stabilizers in $SL_2(\Z)$:
\begin{center}
\begin{tabular}{c r c c}
cusps of $\pm\Gamma^1(5)$ \qquad\qquad & generators& of& stabilizers\\
\hline $\infty$ \qquad\qquad
 & $\gamma$&=&$\displaystyle{\left(\begin{matrix}1&5\\0&1\end{matrix}
                   \right)}$\\
$0$ \qquad\qquad
 & $\delta$&=&$\displaystyle{\left(\begin{matrix}1&0\\-1&1\end{matrix}
                 \right)}$\\
$-2$ \qquad\qquad
 & $A\gamma
A^{-1}$&=&$\displaystyle{\left(\begin{matrix}11&20\\-5&-9\end{matrix}
                 \right)}$\\
$-\frac52$ \qquad\qquad
 & $A\delta
A^{-1}$&=&$\displaystyle{\left(\begin{matrix}11&25\\-4&-9\end{matrix}
                 \right)}$\\
\end{tabular}
\end{center}

\noindent Therefore the group $\Gamma^1(5)$ is generated by
$\gamma$, $\delta$, $A\gamma A^{-1}$, $A\delta A^{-1}$ with the
relation
$$(A\delta
A^{-1})(A\gamma A^{-1})\delta \gamma = I.$$ In particular,
$\Gamma^1(5)$ is actually generated by $A\delta A^{-1}$, $\delta$,
$\gamma$.

For $\Gamma$, we may assume that its cusp at $\infty$ has width 15
so that
$$\Gamma^1(5) = \bigcup_{0 \le j \le 2} \Gamma \gamma^j .$$
\noindent The information on types of special fibers of $\E$ and
the above table give rise to the following information on cusps of
$\Gamma$ and a choice of generators of stabilizers of each cusp:

\begin{center}
\begin{tabular}{c c c}
cusps of $\Gamma$ \qquad\qquad & width \qquad\qquad & generators
of
stabilizers\\
\hline
$\infty$     \qquad\qquad & $15$ \qquad\qquad & $\gamma^3$\\
$-2$         \qquad\qquad & $15$ \qquad\qquad & $A\gamma^3 A^{-1}$\\
$0$          \qquad\qquad & $1 $ \qquad\qquad & $\delta$\\
$5$          \qquad\qquad & $1 $ \qquad\qquad &
$\gamma\delta\gamma^{-1}$\\
$10$         \qquad\qquad & $1 $ \qquad\qquad &
$\gamma^2\delta\gamma^{-2}$\\
$-\frac52$   \qquad\qquad & $1 $ \qquad\qquad & $A\delta A^{-1}$\\
$\frac52$    \qquad\qquad & $1 $ \qquad\qquad & $\gamma A\delta
A^{-1}\gamma^{-1}$\\
$\frac{15}{2}$ \qquad\qquad & $1$ \qquad\qquad & $\gamma^2 A\delta
A^{-1}\gamma^{-2}$\\
\end{tabular}
\end{center}

\noindent This shows that $\Gamma$ is generated by $\gamma^3$,
$\delta$, $A\gamma^3 A^{-1}$, $A\delta A^{-1}$, $\gamma\delta
\gamma^{-1}$, $\gamma A\delta A^{-1}\gamma^{-1}$,
$\gamma^2\delta\gamma^{-2}$, $\gamma^2 A\delta A^{-1}\gamma^{-2}$
with the relation
\[
(A\delta A^{-1})(A\gamma^3 A^{-1}) \delta (\gamma A\delta
A^{-1}\gamma^{-1}) (\gamma\delta\gamma^{-1})(\gamma^2 A\delta
A^{-1}\gamma^{-2}) (\gamma^2\delta\gamma^{-2})\gamma^3 = I.
\]

\begin{rem}\label{r:2.1}
Similar to the above, if we take elements $\gamma^2$, $\delta$,
$A\delta A^{-1}$, $A\gamma^2 A^{-1}$, $\gamma\delta\gamma^{-1}$,
$\gamma A\delta A^{-1}\gamma^{-1}$ as generators with the relation
\[
(A\delta A^{-1})(A\gamma^2 A^{-1})(\delta) (\gamma A\delta A^{-1})
(\gamma\delta\gamma^{-1})(\gamma^2)=I,
\]
then we get the noncongruence subgroup $\G_2$ of $\Gamma^1(5)$ of
index $2$ associated to the elliptic modular surface defined by
\begin{equation}\label{eqn_deg2}
  y^2+(1-t^2)xy-t^2y=x^3-t^2x^2\, .
\end{equation}
\end{rem}

\section{The space of weight 3 cusp forms for $\Gamma$}\label{cusp_forms}

It follows readily from the dimension formula in Shimura
\cite{shimu} that the space $S_3(\Gamma)$ of cusp forms of weight
3 for $\Gamma$ has dimension 2. We shall give a basis of this
space in terms of the weight $3$ Eisenstein series of
$\Gamma^1(5)$ by using Hecke's construction as described in Ogg
\cite{Ogg}. The space of weight $3$ Eisenstein series of
$\Gamma^1(5)$ has dimension $4$, equal to the number of cusps of
$\Gamma^1(5)$. We are only interested in the two Eisenstein series
that vanish at all but only one of the two cusps $\infty$ and
$-2$.

Let $k\ge 3$, and $c, d\in \Z$. The Eisenstein series
\[
G_k(\tau;(c,d);N)=\sum\ ^{'}_{\substack{m\equiv c \pmod N\\n\equiv
d \pmod N}
   } \left(m\tau + n\right)^{-k}
\]
is a weight $k$ modular form for the principal congruence subgroup
$\Gamma(N)$. Moreover, for any $L\in SL_2(\Z)$, we have
\[
G_k(\tau;(c,d);N)|L = G_k(\tau;(c,d)L;N).
\]
This Eisenstein series has the following Fourier expansion:
\begin{equation}\label{e:3.1}
G_k(\tau;(c,d);N)=\sum^\infty_{\l=0} a_{\l} z^{\l}, \qquad\qquad
z=e^{2\pi i\tau/N}
\end{equation}
where
\[
a_0=
\begin{cases}
0 \qquad\qquad\qquad &\text{ if } c\not\equiv 0\,\pmod N,\\
\sum_{n\equiv d\,\pmod N} n^{-k} &\text{ if } c\equiv 0 \,\pmod N,
\end{cases}
\]
and for $\l\ge 1$,
\begin{equation}\label{e:3.2}
a_{\l}=\frac{\left(-2\pi i\right)^k}{N^k \Gamma(k)}
   \sum_{\substack{m\nu =\l\\m\equiv c\,\pmod N}} (\sgn \nu)
   \nu^{k-1} e^{2\pi i\nu d/N}.
\end{equation}
The restricted Eisenstein series is defined by
\[
G^*_k(\tau;(c,d);N)=\sum_{\substack{m\equiv c \pmod N\\ n\equiv d
\pmod N
    \\(m,n)=1}} \left(m\tau +n\right)^{-k},
\]
which is a weight $k$ modular form for $\Gamma(N)$, satisfying
\[
G^*_k(\tau;(c,d);N)|L = G^*_k(\tau;(c,d)L;N)
\]

\noindent for all $L\in SL_2(\Z)$. Let $\mu(n)$ denote the
M\"obius function. Then $G^*_k(\tau;(c,d);N)$ can be expressed in
terms of the series $G_k(\tau; (c,d);N)$:
\begin{equation}\label{e:3.3}
\begin{split}
G^*_k(\tau;(c,d);N)&=\sum^\infty_{a=1}\mu(a) a^{-k}
G_k(\tau;(a'c,a'd);N)\\
     &=\sum_{\substack{(t,N)=1\\t \mod N}}c_t\cdot
G_k(\tau;(ct,dt);N),
\end{split}
\end{equation}

\noindent where $a'$ is chosen such that $a a'\equiv 1 \pmod N$,
and $c_t =\sum_{at\equiv 1 \pmod N, a>0} \mu(a) a^{-k}$. The
Eisenstein series $G^*_k(\tau;(c,d);N)$ has value 1 at the cusp
$-\frac dc$ and 0 at all other cusps.

In the special case $k=3$, $N=5$, for any character $\chi$ of
$\left( \Z/5\Z\right)^*$ we have
\[
\sum_{(t,N)=1}\bar\chi(t)c_t=\sum_t\sum_{a\equiv
t^{-1},\,a>0}\mu(a) a^{-k} =L^{-1}(k,\chi).
\]
Denote by $\chi_3: \left(\Z/5\Z\right)^*\to \C$ the character
given by $\chi_3(2) = i$, $\chi_2=\chi_3^2$, $\chi_4=\chi_3^3$,
and
 $\chi_1$ the trivial character of
$\left(\Z/5\Z\right)^*$. Then the constants $c_t$ can be expressed
via the values of $L$-series:
\[
\begin{split}
&c_1=\frac{1}{4}\left(L^{-1}(3,\chi_1)+L^{-1}(3,\chi_2)+L^{-1}(3,\chi_3
)
    +L^{-1}(3,\chi_4)\right),\\
&c_2=\frac{1}{4}\left(L^{-1}(3,\chi_1)-L^{-1}(3,\chi_2)+iL^{-1}(3,\chi_
3)
    -iL^{-1}(3,\chi_4)\right),\\
&c_3=\frac{1}{4}\left(L^{-1}(3,\chi_1)-L^{-1}(3,\chi_2)-iL^{-1}(3,\chi_
3)
    +iL^{-1}(3,\chi_4)\right),\\
&c_4=\frac{1}{4}\left(L^{-1}(3,\chi_1)+L^{-1}(3,\chi_2)-L^{-1}(3,\chi_3
)
    -L^{-1}(3,\chi_4)\right).
\end{split}
\]
Using the functional equation of the $L$-function $L(s,\chi)$ and
the Bernoulli polynomials, we obtain two explicit $L$-values
\[
\begin{split}
&L(3,\chi_3)=\frac{\tau(\chi_3)}{2i}\left(-\frac{1}{2}\right)
  \left(\frac{2\pi}{5}\right)^3\left(-\frac{1}{3}\right)
  \frac{6}{5}\left(2-i\right),\\
&L(3,\chi_4)=\frac{\tau(\chi_4)}{2i}\left(-\frac{1}{2}\right)
  \left(\frac{2\pi}{5}\right)^3\left(-\frac{1}{3}\right)
  \frac{6}{5}\left(2+i\right),\\
\end{split}
\]
where $\tau(\chi)$ denotes the Gauss sum of the character $\chi$.

As $\Gamma^1(5)=\bigcup^4_{a=0} \Gamma(5)\displaystyle{\left(
    \begin{matrix}1 &0\\a& 1\end{matrix}\right)}$, we let
\begin{equation}\label{e:3.5}
\begin{split}
E_1(\tau)&=\sum_{0 \le a \le 4} G^*_k(\tau;(0,1);5)|\left(\begin{matrix}1 &0\\
                                            a &1\end{matrix}\right)
=\sum_{a}G^*_k(\tau;(a,1);5)\\
&=c_1\sum_a G_k(\tau;(a,1);5)+c_2\sum_a G_k(\tau;(2a,2);5)
  +c_3\sum_a G_k(\tau;(3a,3);5) \\
 \qquad\qquad &+ c_4\sum_a G_k(\tau;(4a,4);5).
\end{split}
\end{equation}
Applying the Fourier expansion of the Eisenstein series
$G_k(\tau;(c,d);5)$ and the equations for $c_i's$, we have the
Fourier expansion of $E_1(\tau)$:
\begin{equation}\label{eqn_E1}
E_1(\tau)=1-\frac{1}{2}\sum^\infty_{\l=1}\left(
2\sum_{\nu|\l,\,\nu>0}
  \nu^2\left(\chi_3(\nu)+\chi_4(\nu)\right) + i \sum_{\nu|\l,\,\nu>0}
  \nu^2\left(\chi_4(\nu)-\chi_3(\nu)\right) \right) q^{\l/5}.
\end{equation}
In particular, the Fourier coefficients of $E_1(\tau)$ are
rational integers. And from the definition of $E_1(\tau)$, it has
value 1 at the cusp $\infty$, and 0 at the other cusps 0, $-2$,
$-\frac52$.

Let
\[
\begin{split}
E_2(\tau)&=E_1(\tau)|A^{-1}\\
    &=c_1\sum_a G_k(\tau;(-2a+1,-5a+2);5)+c_2\sum_a
G_k(\tau;(-4a+2,-10a+4);5)\\
    &+c_3\sum_a G_k(\tau;(-6a+3,-15a+6);5)+c_4\sum_a
G_k(\tau;(-8a+4,-20a+8);5).
\end{split}
\]
Then $E_2(\tau)$ has value $-1$ at the cusp $-2$, and 0 at other
cusps. In the same way we calculate the Fourier expansion of this
Eisenstein series to get
\begin{equation}\label{eqn_E2}
E_2(\tau)=\frac{1}{2}\sum^\infty_{\l=1}\left(
\sum_{\nu|\l,\,\nu>0} \nu^2\left(
  \chi_3(\nu)+\chi_4(\nu)\right)+2i\sum_{\nu|\l,\,\nu>0} \nu^2\left(
\chi_3(\nu)
  -\chi_4(\nu)\right) \right) q^{\l/5}.
\end{equation}

Since $E_1(\tau)$ has values 1, 0, 0, 0 at the cusps $\infty$,
$-2$, $0$, $-\frac52$, respectively, and $E_2(\tau)$ has values 0,
$-1$, 0, 0 at these cusps, both modular forms have no other zero
points. Consider the natural covering map
\[
\Gamma\backslash\CH \to \Gamma^1(5)\backslash \CH,
\] where $\CH$ denotes  the Poincar\'e upper half plane.
It ramifies  only at the two cusps $\infty$ and $-2$, with index
3. Therefore the two functions
\begin{equation}
f_1=\sqrt[3]{E_1^2(\tau) E_2(\tau)} \qquad\qquad \text{ and
}\qquad\qquad f_2=\sqrt[3]{E_1(\tau)E_2^2(\tau)}
\end{equation}
are well-defined entire modular forms of weight 3 for $\Gamma$.
Further, they vanish at every cusp, hence they are cusp forms. It
is clear that $f_1$ and $f_2$  are linearly independent, and thus
form a basis of the space $S_3(\Gamma)$. Choosing a proper cubic
root of one, we may assume that the Fourier coefficients of both
$f_1$ and $f_2$ are rational numbers with denominators involving
only powers of 3.

Since the action of the matrix $A$ interchanges the cusp $\infty$
with the cusp $-2$, it defines an operator on the space
$S_3(\Gamma)$. More precisely, its actions on $f_1$ and $f_2$ are
\[
(f_1)|A=f_2 \qquad\qquad \text{ and }\qquad\qquad (f_2)|A=-f_1.
\]
Thus the operator $A$ on $S_3(\Gamma)$ has eigenforms $f_+
=f_1+if_2$ and $f_- =f_1-if_2$ with eigenvalues $-i$ and $i$,
respectively. The  Fourier expansions  of these two eigenforms are
as follows:

\begin{eqnarray*}
f_+(\tau) &=&  q^{1/15} + iq^{2/15} - \frac {11}{3}q^{4/15}
-i\frac{16}{3}q^{5/15} -\frac {4}{9}q^{7/15} + i\frac
{71}{9}q^{8/15} +\frac{932}{81}q^{10/15}\\
&&+  i\frac{247}{81}q^{11/15}+
\frac{443}{243}q^{13}-i\frac{3832}{243}q^{14/15}-\frac{13151}{729}q^{16/15}+i\frac{9131}{729}q^{17/15}+O(q^{18/15}),
 \\
f_-(\tau) &=& q^{1/15} - iq^{2/15} - \frac {11}{3}q^{4/15}
+i\frac{16}{3}q^{5/15} -\frac {4}{9}q^{7/15} - i\frac
{71}{9}q^{8/15} + \frac{932}{81}q^{10/15}\\ &&-
i\frac{247}{81}q^{11/15}+
\frac{443}{243}q^{13}+i\frac{3832}{243}q^{14/15}-\frac{13151}{729}q^{16/15}-i\frac{9131}{729}q^{17/15}+O(q^{18/15}).
\end{eqnarray*}
This proves the first assertion of Theorem 1.2.

Let $X$ be the modular curve of the noncongruence subgroup
$\Gamma$, that is, $X(\C)=\overline{\Gamma\backslash\CH}$. As seen
in \S \ref{some_elliptic_surface}, it is a projective line over
$\C$. The two cusp forms constructed in section \ref{cusp_forms}
give rise to a Hauptmodul of $X$:
\[
t=\frac{f_1}{f_2} = \sqrt[3]{\frac{E_1}{E_2}}.
\]
Since the Fourier coefficients of $t$ are in $\Q$, the curve $X$
is defined over $\Q$. It is easy to check from the generators and
relation exhibited in \S 3 that the matrix $A$ lies in the
normalizer of the noncongruence subgroup $\Gamma$ in $SL_2(\Z)$.
Therefore $A$ induces a $\Q$-rational involution on the modular
curve $X$, given by
\begin{equation}\label{e:4.1}
A(t)=-\frac{1}{t}.
\end{equation}
Further, $A$ induces an order 4 $\Q$-rational action on the
elliptic surface $\E$.  We first notice that the $j$-function
(\ref{j}) is invariant if we send $t$ to $-1/t$. To see the action
more explicitly, make the following change of variables:
\begin{eqnarray*}
  x&=&t^3X-1/12t^6+1/2t^3-1/12\\
  y&=&t^4Y+1/2t^6X-1/2t^3X-1/24t^9+7/24t^6+5/24t^3+1/24
\end{eqnarray*}
so that the original defining equation (\ref{e:1.1}) becomes
\begin{equation}\label{1515_2}
  Y^2=t(X^3-\frac{1+12t^3+14t^6-12t^9+t^{12}}{48t^6}X
  +\frac{1+18t^3+75t^6+75t^{12}-18t^{15}+t^{18}}{864t^9} ).
\end{equation}
The action of $A$ sends $t$ to $-1/t$, $Y$ to $Y/t$, and $X$ to
$-X$. Hence it is defined over $\Q$ and has order 4.

\bigskip

\section{The $l$-adic representation attached to
$S_3(\Gamma$)}\label{l-adic_repn}

As explained in \S1, given a noncongruence subgroup $\G'$ of
$SL_2(\Z)$ of finite index with the modular curve $X_{\G'}$
defined over $\Q$,  Scholl in  \cite{Sch85} defined a compatible
family of $l$-adic representations $\rho_l$ of
$\Gal(\overline\Q/\Q)$ associated to $S_k(\G')$, $k \ge 3$, from
which he derived the congruence relation Theorem 1.2. When $k =
3$, the representation $\rho_l$ is defined as follows. Choose an
integer $N \ge 3$ such that $\pm\G' \G(N) = SL_2(\Z)$. Denote by
$X(N)$ the compactified modular curve for the principal congruence
subgroup $\Gamma(N)$, and by $X(N)^o$ the part of $X(N)$ with
cusps removed. Let $G(N)=SL(\mu_{N}\times\Z/N)$ , let $f^{\rm
univ}: E^{\rm univ}\to X(N)^{o}$ be the restriction to $X(N)^{o}$
of the universal elliptic curve of $X(N)$, and let $V(N)$ (resp.
$V(N)^o$) be the normalization of the fiber product
$X_{\G'}\times_{X(1)} X(N)$ (resp. $X_{\G'}\times_{X(1)} X(N)^o$).
The finite group scheme $G(N)$ acts on the second factor of
$V(N)$, $E^{\rm univ}$, and the sheaf $\CF_l^{\rm
univ}=R^1f_*^{\rm univ} \Q_l$, respectively. We have the
projection map $\pi'_0:V(N)^{o}\to X(N)^{o}$ and the inclusion map
$i_N: V(N)^{o}\to V(N)$. The representation $\rho_l$ of
$\Gal(\overline\Q/\Q)$ is the action of the Galois group on the
$\Q_l$-space
\[
H^1(V(N)\otimes\overline\Q, (i_N)_{*}\pi_0'^* \CF^{\rm
univ}_l)^{G(N)}.
\]
\smallskip

The reason that an auxiliary modular curve $X(N)$ is involved is
that the curve $X_{\G'}$ does not have a universal elliptic curve,
while $X(N)$ for $N \ge 3$ does. As shown above, the $l$-adic
sheaf comes from this universal elliptic curve. At the end, $G(N)$
invariants are taken to rid the dependence of $X(N)$.
\medskip

When there is an elliptic surface $\E'$ over the modular curve
$X_{\G'}$ with
 $h': \E' \rightarrow X_{\G'}$ tamely ramified along the cusps and elliptic
points, inspired by \cite{Sch88},  we introduce another $l$-adic
representation $\rho_l^*$ of the Galois group of $\Q$ using $\E'$
as follows. Let $X_{\G'}^0$ be the part of $X_{\G'}$ with the
cusps and elliptic points removed. Denote by $i$ the inclusion
from $X_{\G'}^0$ into $X_{\G'}$, and by
\[
h: \E' \rightarrow X_{\G'}^0
\]
the restriction map, which is smooth. For almost all prime $l$, we
obtain a sheaf
\[
\CF_l=R^1h_* \Q_l
\]
on $X_{\G'}^0$. The $i_*$ map then transports it to a sheaf
$i_*\CF_l$ on $X_{\G'}$.  The action of $\rm{Gal}(\overline\Q/\Q)$
on the $\Q_l$-space
\begin{equation}\label{e:4.2}
W_l=H^1(X_{\G'}\otimes\overline\Q,i_*\CF_l)
\end{equation}
\noindent defines an $l$-adic representation, denoted by
$\rho_l^{*}$, of the Galois group of $\Q$.

 The following diagram depicts the relationship of the curves and
surfaces involved. For a scheme $X$ and a non-zero integer $M$, we
use $X[1/M]$ to denote $X \times _{{\rm Spec} \Z} {\rm Spec}
\Z[1/M]$. The integer $M$ below is chosen so that both modular
curves $X(N)$ and $X_{\G'}$ are smooth and proper over the ring
$\Z[1/M]$. The maps $\pi'_0$ and $\pi_{\G'}$ are the natural
projections.
\medskip

$$ \xymatrix{ E^{univ} \ar[d]^{f^{univ}}&&&\E' \ar[d]^{h}\\
X(N)^0  &V^o(N)[1/Ml] \ar[l]_{\pi_0^*}\ar[r]^{i_N} \ar[d]&
V(N)[1/Ml] \ar[d]
\ar[dr]^{\pi_{\G'}}&X_{\G'}^o [1/Ml] \ar[d]^{i} \\
&X_{\G'}\times_{X(1)} X(N)^o  [1/Ml] \ar[r]&X_{\G'}\times_{X(1)}
X(N)  [1/Ml] \ar[r] &X_{\G'}[1/Ml]}$$

\begin{prop}\label{eqrepn}
The two representations $\rho_l^{*}$ and $\rho_l$ are isomorphic
up to a twist by a character $\phi_l$ of $\Gal(\overline\Q/\Q)$ of
order at most 2.
\end{prop}
\begin{proof}
Using the projection $\pi_{\Gamma'}$ from the fiber product
$V(N)[1/Ml]$ to the factor $X_{\G'}[1/Ml]$, we pull back the sheaf
$i_*\CF_l$ on $X_{\G'}[1/Ml]$ to the sheaf $\pi_{\G'}^{*}i_*\CF_l$
on $V(N)[1/Ml]$. It follows from an argument similar to  \S1.3 of
\cite{Sch88} that the sheaf $(i_N)_{*} \pi_0'^* \CF^{\rm univ}_l$
is isomorphic to the sheaf $\pi_{\G'}^{*}i_*\CF_l \otimes \CL$,
where $\CL$ is a rank one sheaf on $V(N)[1/Ml]$ with $\CL^{\otimes
2}$ isomorphic to the constant sheaf $\Q_l$. Consequently,
$H^1(V(N)\otimes\overline\Q,\pi_{\G'}^{*}i_*\CF_l)$ and
$H^1(V(N)\otimes\overline\Q, (i_N)_{*}\pi_0'^* \CF^{\rm univ}_l)$
are isomorphic $\rm{Gal}(\overline\Q/\Q)$ modules up to twisted by
a character $\phi_l$ of $\Gal(\overline\Q/\Q)$ of order at most 2.
Since $G(N)$ acts only on $X(N)$, the $G(N)$-invariant part of
$H^1(V(N)\otimes\overline\Q, \pi_{\G'}^{*}i_*\CF_l)$ is isomorphic
to $W_l$. Therefore the representation $\rho_l^{*}$ on $W_l$ is
isomorphic to $\rho_l$ on $H^1(V(N)\otimes\overline\Q,
(i_N)_{*}\pi_0'^* \CF^{\rm univ}_l)^{G(N)}$ up to twisted by
$\phi_l$.

\end{proof}
\medskip

Apply the above discussion to the case where $\G' = \G$. We shall
show later in \S 6 that our $\rho_l^*$ is in fact isomorphic to
Scholl's representation $\rho_l$. The simpler description of
$\rho_l^*$ allows us to get more information about the
representation, and eventually leading to a finer congruence
result than the one provided by Theorem 1.1.

By Scholl's result in \cite{Sch85}, $\dim_{\Q_l} W_l= h^{2,0}(\E)
+ h^{0,2}(\E)=4$. As remarked at the end of the previous section,
the action of $A$ on $X$ and on $\E$ are both $\Q$-rational, thus
$A$ commutes with the action of ${\rm Gal}(\overline \Q/\Q)$ on
 the space $W_l$. Moreover, the action of $A$ on the
sheaf $\CF_l$ or the representing space $W_l$ has order 4. This
makes $W_l$ a 2-dimensional module over the algebra $K = \Q_l(A)$.

Now fix the prime $l=2$. Then $K= \Q_2(A)$ is a degree two field
extension of $\Q_2$ and $\rho_2^{*}$ is a degree two
representation of ${\rm Gal}(\overline \Q/\Q)$ over the
2-dimensional vector space $W_2$ over $K$. We calculate the
characteristic polynomial $H_p(T)$ (resp. $H'_p(T)$) of
$\rho_2^{*}(\Frob_p)$ over $\Q_2$ (resp. $\Q_2(A) =K =
\Q(A)_{1+A}$) for varying primes $p\ne 2,3$. By Scholl's work
\cite{Sch85}, the characteristic polynomial $H_p(T)$ of
$\rho_2^{*}(\Frob_p)$ can be factored as
\begin{equation}\label{e:4.3}
\begin{split}
H_p(T)&=(T-\a_p)(T-\b_p)(T-p^2/\a_p)(T-p^2/\b_p)\\
      &=T^4-C_1(p)\,T^3+C_2(p)\,T^2-p^2C_1(p)\,T+p^4 \in \Z[T],
\end{split}
\end{equation}
where
\[
\begin{split}
&C_1(p)=\a_p+\b_p+p^2/\a_p+p^2/\b_p=\Tr(\rho_2^{*}(\Frob_p)),\\
&C_2(p)=\frac{1}{2}\left(C_1^2-\Tr(\rho_2{*}(\Frob_p^2))\right)=
\frac{1}{2}\left(\left(\Tr(\rho_2^{*}(\Frob_p))\right)^2-\Tr(\rho_2^{*}(\Frob_{p^2}))\right).
\end{split}
\]

Since
$$
H^j(X_{\G'}\otimes\overline{\F}_p,i_*\CF_l)=0 \quad {\rm for }\,\,
j\neq 1,
$$
 we can easily get from the Lefschetz formula the following trace
formula of $\Tr(\rho_2^{*}(\Frob_q))$ for $q=p^r$:
\[
\Tr(\rho_2^{*}(\Frob_q))=-\sum_{x\in X(\F_q)} \Tr(x),
\]
where
\[
\Tr(x)= \Tr\left((\Frob_q)_x:(i_*\CF_2)_x\right)
\]
\noindent is the trace of $\Frob_q$ restricted to the stalk at $x$
of the sheaf $i_*\CF_2$.

We proceed to compute $\Tr(x)$. Denote by $\E_x$ the fibre of the
elliptic surface $f: \E\to X$ at the point $x \in X(\F_q)$; it is a
curve of genus 0 or 1 depending on whether the discriminant of
$\E_x$ vanishes in the field $\F_q$ or not. Note that there is no
need to treat the cases
 $j=0$ and $j=1728$ separately as in \cite{Sch88} since
 we are computing the traces by using the explicit equation of
the elliptic surface $\E\to X$.
In our case, for $x\in X(\F_q)$ we always have
\[
\Tr(x)=1+q-\#\E_x(\F_q).
\]

\noindent A computer program  yields the following table for
$\Tr_q$, the traces of $\Frob_q$:

\begin{center}\label{table_tr}
\begin{tabular}{|c|c|c|c|c|c|c|c|c|c|}
\hline
$p$         & 5 & 7 & 11 & 13 & 17 & 19 & 23 & 29 & 31 \\
\hline
$\Tr_{p}$ & 0 &10 &0  &$-20$& 0& $-32$& 0  & 0 &$-2$ \\
\hline
$\Tr_{p^2}$&82 &$-146$ &34 &$-476$ &508 &$-932$ &1828 &1564 &$-3842$\\
\hline
\end{tabular}\\
\end{center}
Therefore we obtain the characteristic polynomials $H_p(T)$ of
$\rho_2^{*}(\Frob_p)$, with $p$ the primes between 5 and 31. For
such a prime $p$, the characteristic polynomial $H'_p(T)$ of
$\rho_2^{*}(\Frob_p)$ over $K$ has degree 2, and it has the
property that $H'_p(T)H''_p(T) = H_p(T)$, where  $H''_p(T)$ are
the conjugate of $H'_p(T)$ under the automorphism of $K$ over
$\Q_2$ sending $A$ to $-A$. Since Scholl [sch85] proved that all
roots of $H_p$ are algebraic integers, the first step towards
determining $H'_p$ is to figure out how to separate the four roots
of $H_p$ into two conjugate pairs to form the roots of $H'_p$ and
$H''_p$. It turns out that for the primes from 5 to 31 such
separation is unique except for $p = 13$ and 19. For these primes,
while we cannot choose between $H'_p$ and $H''_p$ without further
work, we do know the determinants of the two-dimensional
representation $\rho_2^{*}$ over $K$ at these primes since the
constant term of $H'_p(T)$ is $\pm p^2$.

Next we prove that the information so far determines the
determinants of $\rho_2^{*}$, which in turn will allow us
determine $H'_{13}$ and $H'_{19}$.

\medskip

\begin{lemma}\label{l:5.1}
If two integral 1-dimensional representations $\sigma_1$ and
$\sigma_2$ of $\Gal(\overline\Q/\Q)$ over the field $\Q_2(A)$,
which are unramified away from 2 and 3, agree on the elements
$\Frob_p$ for $p=5, 7, 11, 17$, then they are equal.
\end{lemma}
\begin{proof}
The images of the 1-dimensional representations are in
$\left(\Z_2[A] \right)^*=\left(\Z(A)_{\p}\right)^*$, where
$\p=(1+A)$ is the maximal ideal of the local ring $\Z_2(A)$. Note
that $\left(\Z(A)_{\p}\right)^*=\la A\ra\times
\left(1+\p^3\right)$. Let $\log$ denote the $(1+A)$-adic logarithm
on $\left(\Z_2[A]\right)^*$. More precisely, it has kernel the
group $\la A\ra$ of roots of unity in $\Q_2(A)$, and it maps $1 +
x \in 1+\p^3$
 to $\sum^\infty_{n=1} \frac{(-1)^{n-1}}{n} x^n \in \mathbb{Z}_2[A]$.
As such, $\log$ gives an isomorphism between the multiplicative
group $1+\p^3$ and the additive group $\p^3$. Consider
\[
\psi= \log\circ \sigma_1 -\log\circ\sigma_2,
\]
which is a homomorphism from $\Gal(\overline\Q/\Q)$ to $\p^3$.

If $\psi\ne 0$, then
\[
n_0={\rm min}\left\{ {\rm ord}_{\p}\left(\psi(\tau) \right):
\tau\in \Gal(\overline\Q/\Q)\right\}
\]
is finite. Then
$$\bar\psi :=\frac{1}{(1+A)^{n_0+1}} \psi \pmod \p$$
\noindent is a continuous surjective homomorphism from
$\Gal(\overline\Q/\Q)$ to $\F_2$ which is trivial at the Frobenius
elements at primes $ p = 5, 7, 11, 17$ by assumption. This
representation of $\Gal(\overline\Q/\Q)$ factors through a
quadratic extension of $\Q$ unramified outside $2$ and $3$. Such
fields are extensions of $\Q$ by adjoining square roots of $2, 3,
6, -1, -2, -3, -6$, respectively. It is easy to check that the
prime $ p = 5, 5, 11, 7, 5, 5, 17$ is inert in the respective
field, and thus $\bar\psi$ at such $\Frob_p$ would be nontrivial,
a contradiction. Therefore $\psi = 0$, in other words, the image
of the representation $\sigma :=\sigma_1(\sigma_2)^{-1}$ is a
subgroup of $\la A\ra$, a cyclic group of order 4. Hence we
consider all Galois extensions of $\Q$ with group equal to a
subgroup of a cyclic group of order 4,
 unramified away from $2$ and $3$, and in which $p = 5, 7, 11, 17$
split completely. If a nontrivial such extension exists, then it
contains a quadratic subextension unramified outside 2 and 3, and
in which $p = 5, 7, 11, 17$ split completely. As shown above, this
is impossible. Therefore the image of $\sigma$ can only be
$\{1\}$, in other words, $\sigma_1$ and $\sigma_2$ are equal.
\end{proof}

Denote by $\chi_{-3}$ the quadratic character attached to the
field $\Q(\sqrt {-3})$, that is, $\chi_{-3}(x)$ is equal to the
Legendre symbol $\left(\frac {-3}{x}\right)$. The 1-dimensional
representation $\sigma$ of $\Gal(\overline\Q/\Q)$ over $K$ given
by $\tilde \chi_{-3}(\Frob_p) = \chi_{-3}(p)p^2$ for primes $p \ne
2, 3$ agrees with $\det(\rho_2^{*})$ at $\Frob_p$ for $p = 5, 7,
11, 17 $ by checking the constant term of $H'_p(T)$, and hence we
conclude from Lemma \ref{l:5.1} that the two degree one
representations agree.

\medskip
\begin{cor}\label{c:5.1} Let $\tilde \chi_{-3}$ be as above.  We have
$\det(\rho_2^{*}) = \tilde \chi_{-3}$.
\end{cor}

\medskip

In particular, we know that the constant term of $H'_{13}$ (resp.
$H'_{19}$) is $(13)^2$ (resp. $(19)^2$). This information enables
us to determine $H'_{13}(T)$ and $H'_{19}(T)$. We record the
result so far in the following proposition.
\medskip

\begin{prop} The characteristic polynomials $H_p(T)$ and $H'_p(T)$  of
$\rho_2^{*}(\Frob_p)$ over $\Q_2$ and $K = \Q_2(A) = \Q(A)_{1+A}$,
respectively,  for primes $5 \le p \le 31$ are as follows.
\begin{center}
\begin{tabular}{|c|c|c|}
\hline
$p$   & $H_p(T)$    &   $H'_p(T)$ \\
\hline \hline
5  & $T^4-41\,T^2+625$ & $T^2\pm 3A\,T-25$\\
\hline
7  & $T^4-10\,T^3+123\,T^2-490\,T+7^4$  & $T^2-5\,T+7^2$\\
\hline
11 & $T^4-17\,T^2+11^4$  & $T^2 \pm 15A\,T-11^2$\\
\hline
13 & $T^4+20\,T^3+438\,T^2+20\cdot 13^2\,T+13^4$ & $T^2+10\,T+13^2$\\
\hline
17 & $T^4-254\,T^2+17^4$ & $T^2\pm 18A\,T-17^2$ \\
\hline
19 & $T^4+32\,T^3+978\,T^2+32\cdot 19^2\,T+19^4$ & $T^2+16\,T+19^2$\\
\hline
23 & $T^4-914\,T^2+23^4$  & $T^2 \pm 12A\,T-23^2$\\
\hline
29 & $T^4-782\,T^2+29^4$  & $T^2 \pm 30A\,T-29^2$\\
\hline
31 & $T^4+2\,T^3+1923\,T^2+2\cdot 31^2\,T+31^4$ & $T^2+T+31^2$\\
\hline
\end{tabular}

(Table 1)
\end{center}
In case $p \equiv 2 \pmod 3$, the coefficient of $T$ in $H'_p(T)$
is determined up to sign.
\end{prop}
\bigskip

\section{Comparison with the representation $\tilde\rho_2$ attached to certain
cusp forms in $S_3(\Gamma_1(27))$} \label{s:5}

In the space of weight 3 level 27 cusp forms $S_3(\Gamma_1(27))$,
we find from William A. Stein's website \cite{stein} two Hecke
eigenforms $g_a$ whose $q$ expansion ($q = e^{2\pi iz}$) to order
31 are as follows:

\begin{eqnarray*}g_a&=&q + aq^2 - 5q^4 - aq^5 + 5q^7 - aq^8 + 9q^{10} - 5aq^{11}
- 10q^{13} + 5aq^{14} - 11q^{16} + 6aq^{17}\\&& - 16q^{19} +
5aq^{20} + 45q^{22} - 4aq^{23} + 16q^{25} - 10aq^{26} - 25q^{28} +
10aq^{29} - q^{31} + O(q^{32}),  \end{eqnarray*} where $a$ is a
root of $x^2+9$. The  character of this modular form is
$\chi_{-3}$. Denote by $\tilde\rho_2$ the 4-dimensional $2$-adic
representation of the Galois group of $\Q$ attached to $g_a +
g_{-a}$, established by Deligne [Del]. The  Atkin-Lehner operator
$H_{27}$ acts on the curve $X_1(27)$ as an involution, and it is
$\Q$-rational. Further, it induces an action of order 4 on the
representation space of $\tilde\rho_2$ so that $\tilde\rho_2$ may
be regarded as a 2-dimensional representation over $\Q_2(i) =
\Q(i)_{1+i}$. In view of Corollary 5.1, we have
\medskip

\begin{cor}\label{c:6.1}
$\det(\rho_2^{*})=\det(\tilde\rho_2)$.
\end{cor}
\medskip

Fix an isomorphism from $K$ to $\Q_2(i) = \Q(i)_{1+i}$ such that
the characteristic polynomial of $\rho_2^{*}(\Frob_5)$ agrees with
that of $\tilde\rho_2(\Frob_5)$.

 Denote by $\rho'_2$ the representation $\rho_2*$ viewed over $\Q(i)_{1+i}$.
 Our goal in this section is to show that $\rho'_2$ and $\tilde\rho_2$ are
isomorphic.

To compare two representations from $\Gal(\overline\Q/\Q)$ to
$\GL_2(\Z[i]_{1+i})$, we will apply Serre's method \cite{Serre}.

\begin{theorem}[Serre]\label{t:6.1}
Let $\rho_1$ and $\rho_2$ be representations of
$\Gal(\overline\Q/\Q)$ to $\GL_2(\Z[i]_{1+i})$. Assume they
satisfy the following two conditions:
\begin{itemize}
\item[(1)] $\det(\rho_1)=\det(\rho_2)$; \item[(2)] the two
homomorphisms from $\Gal(\overline\Q/\Q)$ to
    $\GL_2(\F_2)$, obtained from the reductions of $\rho_1$ and
    $\rho_2$ modulo $1+i$, are surjective and equal.
\end{itemize}
If $\rho_1$ and $\rho_2$ are not isomorphic, then there exists a
pair $(\widetilde G, t)$, where $\widetilde G$ is a quotient of
the Galois group $\Gal(\overline\Q/\Q)$ isomorphic to either
$S_4\times \{\pm 1\}$, or $S_4$, or $S_3\times \{\pm 1\}$, and the
map $t: \widetilde G \to \F_2$ has value 0 on the elements of
$\widetilde G$ of order $\le 3$, and 1 on the other elements.
\end{theorem}

\medskip
This result, explained in detail in \cite{Serre} and in a letter
from Serre to Tate, is a specialization of a general idea to
determine the ``deviation" of two non-isomorphic representations
over a local field. When the residue field of the local field is
small, it gives a feasible method to determine an $l$-adic
representation by checking the traces of Frobenii at a small
number of primes. Originated from Faltings work \cite{Fal}, this
idea was made effective by Serre \cite{Serre}. In \cite{Liv}
Livn\'e described how to use Serre's method to prove two
representations with even trace to be isomorphic. We sketch Serre's proof below.

Suppose two representations
 $\rho_1$ and $\rho_2$ are not isomorphic. Then the traces
$\Tr(\rho_1)$ and $\Tr(\rho_2)$ are not identical on  $\mG :=
\Gal(\overline\Q/\Q)$. Write $\pi$ for the uniformizer $1+i$ of
the local field $\Q_2(i)$ for brevity. There exists a highest
power $\pi^n$ for some integer $n$ ($\ge 1$ by condition (2)) such
that
\[
\Tr(\rho_1(s))\equiv \Tr(\rho_2(s)) \mod \pi^n \quad\text{ for
all} ~s\in \Gal(\overline\Q/\Q).
\]
This yields a non-constant (and hence surjective) map
\[
\begin{array}{llrll}
&t: &\mG &\to & \F_2\\
&   &s          &\mapsto
&\left(\Tr(\rho_2(s))-\Tr(\rho_1(s))\right)/\pi^n \mod \pi,
\end{array}
\]
which records the difference between the two representations
$\rho_1$ and $\rho_2$.  Next one seeks to pass this information to
a manageable finite quotient $\tG$ of the Galois group $\mG$ so that
the pair $(\tG, t)$ measures the difference of 
$\rho_1$ and $\rho_2$.  

\smallskip

The three cases of ``deviation" listed in Theorem \ref{t:6.1} were
derived by explicitly computing possible $\tG$ as
follows. 
The hypothesis $\rho_1\equiv \rho_2 \mod \pi^n$ implies
that we may write, for $s \in \mG$,
\[
\rho_2(s) = (1+ \pi^n \,a(s))\,\rho_1(s) \quad\text{ with }
~a(s)\in M_2(\O_K),
\]
showing $t(s) = \Tr(a(s)\rho_1(s)) \mod \pi$. Hence it suffices to
find a quotient $\tG$ capturing $a(s) \mod \pi$
and $\rho_1(s) \mod \pi$ for all $s$ in $G$.  
The relation $\rho_2(s_1s_2) = \rho_2(s_1)\rho_2(s_2)$ yields
$a(s_1s_2) \equiv a(s_1) + \rho_1(s_1)a(s_2)\rho_1(s_1)^{-1} \mod
\pi$ for all $s_1, s_2$ in $\mG$. In other words,
 the map $ s \mapsto a(s) \mod \pi$ from $\mG$ to $M_2(\F_2)$ is a
1-cocycle under the adjoint action of the Galois group on
$M_2(\F_2)$ through $\rho_1$ modulo $\pi$. The map $\theta: s \mapsto
(a(s) \mod \pi, ~\rho_1(s) \mod \pi)$ is a homomorphism
from $\mG$ to the semi-direct product
$M_2(\F_2)\rtimes \GL_2(\F_2)$,
where the group law is $(a, g)\cdot(b, h) = (a + gbg^{-1}, gh)$.
The desired group $\tG$ is isomorphic to the image of $\theta$. It remains to figure
out the possible group structure of $\tG$.  
The projection of $\tG$ to $\GL_2(\F_2) \cong S_3$ is surjective by condition (2). Condition (1)
implies that the trace of $a(s)$ mod $\pi$ is zero. Therefore the
projection of $\tG$ to $M_2(\F_2)$ is a subgroup of the trace zero elements in $M_2(\F_2)$,
which is generated by $\left (
\begin{array}{cc}1&0\\0&1
\end{array} \right ), \left ( \begin{array}{cc}0&1\\0&0
\end{array} \right )$, and $\left ( \begin{array}{cc}0&0\\1&0
\end{array} \right )$. Put together, one finds three possibilities for $\tG$:
they are $\la \left ( \begin{array}{cc}1&0\\0&1
\end{array} \right ), \left ( \begin{array}{cc}0&1\\0&0
\end{array} \right ), \left ( \begin{array}{cc}0&0\\1&0
\end{array} \right )\ra \rtimes \GL_2(\F_2)$, $\la\left (
\begin{array}{cc}1&1\\0&1
\end{array} \right ), \left ( \begin{array}{cc}1&0\\1&1
\end{array} \right )\ra \rtimes \GL_2(\F_2)$, and
$\la \left ( \begin{array}{cc}1&0\\0&1
\end{array} \right )\ra \rtimes \GL_2(\F_2)$,
 corresponding to the three cases given in
the theorem.

\medskip

It follows from Serre's Theorem that two nonisomorphic
representations
 must have different traces at elements
 in $\widetilde G$ of order at least 4. Therefore, to show that two
representations $\rho'_2$ and $\tilde\rho_2$ are isomorphic (and
necessarily ramify at the same places), our strategy is to search
for Galois extensions of $\Q$ with Galois groups isomorphic to
those listed in  Theorem \ref{t:6.1} and unramified where the
representations are unramified. In each of such extensions, if we
can find a Frobenius element of order $\ge 4$ in the Galois group
at which the two representations have the same trace, then they
must be isomorphic.
\medskip

Corollary 6.1 shows that condition (1) holds. We now proceed to
prove that both representations satisfy condition (2) as well.
Note that the residue field of $\Q(i)$ at $1+i$ is $\F_2$, hence
the reductions of $\rho'_2$ and $\tilde\rho_2$ mod $1+i$ yield two
representations of the Galois group of $\Q$ to $GL_2(\F_2)$.
 While we do not know the characteristic polynomial of the Frobenius at almost
all primes $p \ge 5$ for the representation $\rho'_2$, we do know
them modulo $1 + i$ for $ 5 \le p \le 31$ from Table 1, and it is
easy to check that they agree with those from representation
$\tilde\rho_2$ mod $1+i$. Thus condition (2) for our two
representations will follow from

\medskip

\begin{lemma}\label{l:6.1} There is only one representation $\rho$ from
$\Gal(\overline\Q/\Q)$ to $\GL_2(\F_2)$, unramified outside 2 and
3, such that the characteristic polynomial of $\rho(\Frob_p)$ is
equal to $H'_p$ mod $(1+A)$ from Table 1 for primes $p = 5, 7,
13$. Further $\rho$ is surjective.
\end{lemma}

\begin{proof} The existence is obvious. We prove the uniqueness. Note that
$\GL_2(\F_2)$ is isomorphic to the symmetric group on three
letters $S_3$, which is generated by an element of order 2 and an
element of order 3. Denote by $\sigma$ the sign homomorphism from
$S_3$ to $\{\pm 1 \}$, and by $\varepsilon$ the composition
$\sigma \circ \rho$. Then $\varepsilon$ is a character of the
Galois group of $\Q$ of order at most two and it is unramified
outside 2 and 3.

We see from Table 1 that the characteristic polynomials of
$\rho(\Frob_5)$  and $\rho(\Frob_7)$ are $T^2 + T + 1$, hence
$\rho(\Frob_5)$ and $\rho(\Frob_7)$ both have order 3. The
characteristic polynomial of $\rho(\Frob_{13})$ is $T^2 + 1$,
hence $\rho(\Frob_{13})$ has order 1 or 2. In particular,
$\varepsilon(\Frob_p) = 1$ for $p = 5, 7$. If $\varepsilon$ is
nontrivial, then it arises from a quadratic extension of $\Q$
unramified outside 2 and 3. Such extensions are $\Q(\sqrt d)$ with
$d = 2, 3, 6, -1, -2, -3, -6$. Since 5 is inert in $\Q(\sqrt d)$
with $d = 2, 3, -2, -3$ and 7 is inert in $\Q(\sqrt d)$ with $d =
6, -1$, this leaves $\varepsilon = \left( \frac {-6}{} \right)$ or
$\varepsilon = 1$ as the only possibilities.

Assume  $\varepsilon = 1$. Then $\rho$ factors through
$\Gal(\overline\Q/\Q)\twoheadrightarrow C_3\subset \GL_2(\F_2)$.
But the unique $C_3$ extension of $\Q$ unramified outside of
$\{2,3\}$ is $\Q(\zeta_9+\zeta_9^{-1})$, in which primes $p\equiv
\pm 1 \mod 9$ split completely. Hence
\[
{\rm ord}\left(\rho(\Frob_p)\right)=
\begin{cases}
1 &\qquad\qquad\qquad\text{ if }\quad p\equiv \pm 1 \mod 9,\\
3 &\qquad\qquad\qquad\text{ if }\quad p\not\equiv\pm 1\mod 9.
\end{cases}
\]

\noindent This contradicts the fact that $\rho(\Frob_{13})$ has
order at most 2.

Therefore $\varepsilon = \left( \frac {-6}{} \right)$ and $\rho$
is surjective.

Finally, we know that ${\rm Ker}(\rho)=\Gal(\overline \Q/K_6)$ for
some $S_3$-Galois extension $K_6$ of $\Q$ containing the quadratic
field $\Q(\sqrt{-6})$ and with discriminant of type $\pm
2^{\a}3^{\b}$. Such a field $K_6$ is unique by applying the class
field theory to the extension $K_6$ over $\Q(\sqrt{-6})$. And it
is given as the splitting field of a cubic polynomial over $\Q$:
\[
K_6={\rm Split}(x^3+3x-2)
\]
\noindent from the table by H. Cohen in [Coh]. The uniqueness of
$\rho$ then follows from the uniqueness of $K_6$ and the
uniqueness of degree 2 irreducible representation over $\F_2$ of
$S_3$.
\end{proof}

\medskip

Now we are ready to apply Serre's theorem to our representations
$\rho'_2$ and $\tilde\rho_2$. It follows from the above proof that
the fixed field of the possible deviation group $\widetilde G$
contains the field $K_6$. We start by finding all quartic fields
$M$ with $\Gal(M/\Q)=S_4$, which contain $K_6$ and are unramified
outside  $2$ and $3$. There are three such fields. Listed below
are their defining equations, discriminants, and certain primes
$p$ such that $\Frob_p$ is of order $4$ in the group $\Gal(L/\Q)$.

\begin{center}
\begin{tabular}{c c c}
defining equation     &  discriminant & $p$ with order $4$ Frobenius\\
\hline
$x^4-4x-3=0$  &  $(-216)\cdot 8^2=-2^9\cdot 3^3$ & 13, 17, 19, 23 \\
$x^4-8x+6=0$  &  $(-216)\cdot 16^2 = -2^13\cdot 3^3$ & 13, 17 \\
$x^4-12x^2-16x+12=0$ & $(-216)\cdot 16^2=-2^13\cdot 3^3$ & 19, 23\\
\end{tabular}
\end{center}

As $H'_p(T)$ agrees with the characteristic polynomial of
$\tilde\rho_2(\Frob_p)$ for primes $p = 13, 19$, the value of $t$
at such $\Frob_p$ is zero. Hence we may rule out two possible
deviation groups $\widetilde G = S_4\times \{\pm 1\}$ and
$\widetilde G =S_4$. For the case $\widetilde G=S_3\times \{\pm
1\}$, elements of interest are those primes $p \ne 2, 3$ such that

 (6.i) $\Frob_p$ has order 3 in $S_3$, which is equivalent to the
trace of $\Frob_p$ being odd under both representations;  and

(6.ii) There is a quadratic extension $\Q(\sqrt d)$ of $\Q$,
unramified outside 2 and 3, in which $p$ is inert.

\noindent (Consequently, $\Frob_p$ in $\widetilde G$ has order 6.)
For the second statement, we consider $\Q(\sqrt d)$ with $d = 2,
3, 6, -1, -2, -3$ since $\Q(\sqrt {-6})$ is contained in $K_6$ and
$\Gal(K_6/\Q)$ is $S_3$.  As $5$ is inert in $\Q(\sqrt d)$ with $d
= 2, 3, -2, -3$ and $7$ is inert in $\Q(\sqrt d)$ with $d = -1,
6$, and the trace of $\Frob_p$ is odd at $p = 5, 7$ under both
representations, we may take $p$ to be 5 or 7. On the other hand,
$\Frob_p$ has the same trace under $\rho'_2$ and $\tilde\rho_2$
for $p = 5, 7$. Hence the last case of $\widetilde G$ is also
eliminated, and $\rho'_2$ and $\tilde\rho_2$ are isomorphic. We
record this in

\medskip

\begin{theorem}\label{t:6.2}
  The two representations $\rho_2^{*}$ and  $\tilde\rho_2$  are isomorphic.
\end{theorem}

\medskip
It is worth pointing out that this theorem uses the information of
$H_p$ for primes $5 \le p \le 19$ only. Further, the $H'_p(T)$'s
for $p \equiv 1 \pmod 3$ in Table 1 are uniquely determined, given
by the characteristic polynomial of $\tilde\rho_2(\Frob_p)$ as in
the list below.

\medskip
\begin{cor}The characteristic polynomials $H_p(T)$ and $H'_p(T)$  of
$\rho_2^{*}(\Frob_p)$ over $\Q_2$ and $K = \Q_2(A) = \Q(A)_{1+A}$,
respectively,  for primes $5 \le p \le 31$ are as follows.
\begin{center}
\begin{tabular}{|c|c|c|}
\hline
$p$   & $H_p(T)$    &   $H'_p(T)$ \\
\hline \hline
5  & $T^4-41\,T^2+625$ & $T^2 - 3A\,T-25$\\
\hline
7  & $T^4-10\,T^3+123\,T^2-490\,T+7^4$  & $T^2-5\,T+7^2$\\
\hline
11 & $T^4-17\,T^2+11^4$  & $T^2 - 15A\,T-11^2$\\
\hline
13 & $T^4+20\,T^3+438\,T^2+20\cdot 13^2\,T+13^4$ & $T^2+10\,T+13^2$\\
\hline
17 & $T^4-254\,T^2+17^4$ & $T^2 + 18A\,T-17^2$ \\
\hline
19 & $T^4+32\,T^3+978\,T^2+32\cdot 19^2\,T+19^4$ & $T^2+16\,T+19^2$\\
\hline
23 & $T^4-914\,T^2+23^4$  & $T^2 - 12A\,T-23^2$\\
\hline
29 & $T^4-782\,T^2+29^4$  & $T^2 + 30A\,T-29^2$\\
\hline
31 & $T^4+2\,T^3+1923\,T^2+2\cdot 31^2\,T+31^4$ & $T^2+T+31^2$\\
\hline
\end{tabular}

(Table 2)
\end{center}
\end{cor}
\medskip
Finally we prove
\medskip

\begin{theorem}\label{t:6.3} The representations $\rho_2$, $\rho_2^{*}$, and
$\tilde\rho_2$ are isomorphic to each other.
\end{theorem}
\begin{proof} We know from Theorem 6.2 that $\rho_2^{*}$ and $\tilde\rho_2$
are isomorphic. Further, by Proposition 5.1, $\rho_2$ is
isomorphic to $\rho_2^{*} \otimes \phi_2$ for some character
$\phi_2$ of order at most 2. Therefore it remains to determine
$\phi_2$. Since $\rho_2$ and $\rho_2^{*}$ are unramified outside 2
and 3, the character $\phi_2$, if nontrivial, is associated to a
quadratic field $\Q(\sqrt d)$ with $d = 2, 3, 6, -1, -2, -3, -6$.
Write $\phi(p)$ for $\phi_2(\Frob_p)$ for brevity. For odd primes
$p \ge 5$, the characteristic polynomial of $\rho_2(\Frob_p)$ is
$$\tilde H_p(T) := T^4 -C_1(p)\phi(p)T^3 + C_2T^2 - p^2C_1(p)\phi(p)T + p^4,$$
\noindent where
$$ H_p(T) = T^4 -C_1(p)T^3 + C_2T^2 - p^2C_1(p)T + p^4$$
\noindent is the characteristic polynomial of
$\rho_2^{*}(\Frob_p)$. By Theorem 1.1, the cusp form
$$f_+(\tau) = \sum_{n \ge 1}a(n)q^{n/15} = q^{1/15} + iq^{2/15} - \frac
{11}{3}q^{4/15} -i\frac{16}{3}q^{5/15} -\frac {4}{9}q^{7/15} +
i\frac {71}{9}q^{8/15} + \frac{932}{81}q^{10/15} + \cdots ,$$
\noindent in $S_3(\G)$ satisfies the congruence relation
$$\ord_p((a(np^2) - C_1(p)\phi(p)a(np) + C_2(p)a(n) -
p^2C_1(p)\phi(p)a(n/p) + p^4a(n/p^2)) \ge 2(\ord_pn + 1)$$
\noindent for all $n \ge 1$ and $p \ge 5$. Applying this
congruence relation to $n = 1$, $p = 7, 13$ and using the explicit
values of $C_1(p), C_2(p)$ from Proposition 5.2 as well as the
known Fourier coefficients of $f_+$ from \S 4, we find that
$\phi(7) = \phi(13) = 1$. Therefore either $\phi_2$ is trivial or
$\phi_2$ is $\chi_{-3}$, the quadratic character attached to the
field $\Q(\sqrt{-3})$. On the other hand, the two newforms
$g_{3i}$ and $g_{-3i}$ are twist of each other by $\chi_{-3}$,
which in turn implies that the representation $\tilde\rho_2$ is
invariant under twisting by $\chi_{-3}$, and hence so is
$\rho_2^{*}$. Therefore in both cases of $\phi_2$ we have $\rho_2$
isomorphic to $\rho_2^{*}$.
\end{proof}

\section{The Atkin-Swinnerton-Dyer congruence relations}\label{s:7}

 Let $\xi$ be the map on the elliptic surface $\E$ (given by
 (\ref{e:1.1})) sending the base parameter $t$ to $\omega ^2 t$,
 where $\omega = e^{2\pi i/3}$ is a primitive cubic root of one.
It induces an action on the weight 3 cusp forms of $\G$ via
 $\xi(f_1)=\omega f_1, \xi(f_2) = \omega^2 f_2$.

Fix a prime $p \ne 2, 3$. Following the notation in \cite{Sch85},
denote by $F$ the canonical endomorphism of $L_k(X,{\Z}_p)$ coming
from an $F$-crystal with logarithmic singularities. In our case,
$k=1$ so that $k+2=3$ is the weight. Recall that $L_k(X, {\Z}_p)$
is the direct sum of the module $S_{3}(X, {\Z}_p)$ of weight 3
cusp forms with Fourier coefficients in ${\Z}_p$ with its dual
$S_{3}(X, {\Z}_p)^{\vee}$. The basis $\{f_1, f_2\}$ of $S_{3}(X,
{\Z}_p)$ gives rise to a dual basis $\{f^{\vee}_1, f^{\vee}_2\}$
of $S_{3}(X, {\Z}_p)^{\vee}$. Denote the space  $L_k(X,{\Z}_p)$ by
$V$ for brevity. The two operators $A$ and $\xi$ act on $V$ as
follows:

\begin{equation}\label{e1:supp}
\begin{array}{lllll}
&A(f_1) = f_2, &A(f_2) = -f_1, &A(f^{\vee}_1) = f^{\vee}_2,
&A(f^{\vee}_2) = -f^{\vee}_1; \\
&\xi(f_1)=\omega f_1, &\xi(f_2) = \omega^2 f_2, &\xi(f^{\vee}_1) =
\omega^2 f^{\vee}_1, &\xi(f^{\vee}_2) = \omega f^{\vee}_2.\\
\end{array}
\end{equation}
Here the actions of $A$ and $\xi$ on cusp forms were given before,
and their actions on the dual space are obtained by applying the
dual action of an operator $T$:
\[
T(h^{\vee})(v) = h^{\vee}(T^{-1}\cdot v), \quad\text{ for
$h^{\vee}\in S_3(X, \Z_p)^{\vee}$, and $v\in S_3(X, \Z_p)$}.
\]

Denote by $W$ the representation space of the 4-dimensional
$2$-adic representation $\rho_2$ of $\Gal(\overline\Q/\Q)$
described at the beginning of \S 5. The operators $A$, $\xi$ and
$F_p$, the Frobenius at $p$,
 act on $W$.

To prove the Atkin-Swinnerton-Dyer congruence relations at $p$, we
shall compare the $2$-adic theory and $p$-adic theory. We begin by
summarizing the actions of the operators involved.

\begin{prop}\label{p1:supp}
The operators $F_p, F, A, \xi$ satisfy the following relations on
their respective representation spaces.
\begin{itemize}
\item[i.] $AF=FA$, $F_p A = A F_p$;

\item[ii.] $A \xi = \xi^2 A$;

\item[iii.] If $p\equiv 1 \mod 3$, then $\xi F = F \xi$, $\xi F_p
= F_p \xi$; If $p\equiv 2 \mod 3$, then $\xi F = F \xi^2$, $\xi
F_p = F_p \xi^2$, hence in this case, $\xi (AF) = (AF) \xi$, and
$\xi (A\,F_p) = (F_p\,A) \xi$;

\item[iv.] $A^4=1$, $\xi^3 = 1$.
\end{itemize}
\end{prop}

\begin{proof} The assertion (i) follows from the fact that $f_1$ and $f_2$
have rational Fourier coefficients; the assertions (ii) and (iv)
follow from (\ref {e1:supp}). We need only prove (iii). The action
of $\xi$ on the base curve is given by
\[
\xi(t) = \omega^2 \cdot t,
\]
and $\xi$ acts trivially on the elliptic curve over $X$. If $p
\equiv 1 \mod 3$, then
\begin{equation}\label{e1.5:supp}
F_p \, \xi (t)=F_p(\omega^2\cdot t) = (\omega^2)^p\cdot t^p =
\omega^2\cdot t^p = \xi\, F_p (t).
\end{equation}
Therefore on the $2$-adic representation space $W$, the action of
$F_p$ commutes with $\xi$ by functoriality. On the $p$-adic
representation space $V$, the action of $F$ is semi-linear, hence
\[
F\, \xi (t) = F(\omega^2\cdot t) = (\omega^2)^p \cdot
F(t)=\xi\,F(t),
\]
and hence the commutativity of $F$ and $\xi$. In the case $p\equiv
2 \mod 3$, we check the commutativity in the same way as in
equation (\ref{e1.5:supp}):
\[
\begin{array}{ll}
&F_p\,\xi^2(t)=F_p(\omega\cdot t)=(\omega)^p\, F_p(t) =
(\omega)^2\, F_p(t) =
\xi\,F_p(t),\\
&F\,\xi^2(t)=F(\omega\cdot t)=(\omega)^p\, F(t) = (\omega)^2\,
F(t) = \xi\,F(t),
\end{array}
\]
 as described in (iii).

\end{proof}

Denote by $char(U, T)$ the characteristic polynomial of an
operator $T$ on the space $U$. The representation spaces $W$ and
$V$ decompose as direct sums of eigenspaces of the two operators
$A$ and $\xi$, respectively, with eigenvalues appearing as
subindexes:
\begin{subequations}\label{e2:supp}
\begin{align}
&W = W_{i}\oplus W_{-i}, & V = V_{i} &\oplus V_{-i}; \\
&W = W_{\omega} \oplus W_{\omega^2}, &V = V_{\omega} &\oplus
V_{\omega^2}.
\end{align}
\end{subequations}

We have shown in \S 6 that $\rho_2$ is isomorphic to the $2$-adic
representation $\tilde\rho_2$ acting on the 4-dimensional space
$\widetilde W$ attached to the newforms $g_+$ and $g_{-}$ of
weight 3 and level 27. The operator
\[
H_{27} = \left( \begin{matrix}
             0 & -1\\
             27 & 0
             \end{matrix}
            \right)
\]
on $\widetilde W$ plays the same role as $A$ on $W$. Hence the
space $\widetilde W$ decomposes as the sum of two eigenspaces of
$H_{27}$ :

\[
\widetilde W = \widetilde W_{i} \oplus \widetilde W_{-i}.
\]
Denote by $Frob_p$ the action of the Frobenius at $p$ on the space
$\widetilde W$. The isomorphism between $\rho_2$ and $\tilde
\rho_2$ yields immediately

\[
 char(W_{\pm i}, F_p) = char(\widetilde W_{\pm i}, Frob_p).
\]

On the other hand, it was shown in \cite{Sch85} that

\[
char(V,F) = char(W, F_p).
\]
Consequently we have
\[
 char(V_{i}, F)char(V_{-i}, F) = char(\widetilde W_{i}, Frob_p)
char(\widetilde W_{-i}, Frob_p) = char({W_{i}}, F_p)
char({W_{-i}}, F_p).
\]

\medskip
We proceed to prove a refinement of this relation.

\begin{theorem}\label{t2:supp}
There hold
\begin{align}\label{e3:supp}
&char(V_{i}, F) = char(W_{i}, F_p) = char(\widetilde W_{i}, Frob_p), \\
&char(V_{-i}, F) = char(W_{-i}, F_p) = char(\widetilde W_{-i},
Frob_p).
\end{align}
\end{theorem}

\begin{proof}In view of the above analysis, it suffices to prove
$char(V_{i}, F) = char(W_{i}, F_p)$ under the assumption that not
all eigenvalues of $char(V, F)$ are equal. The eigenspaces of $V$
occurred in (\ref{e2:supp}) have the following explicit bases:
\begin{subequations}\label{e4:supp}
\begin{align}
&V_i=\la f_1-i f_2, f^{\vee}_1-i f^{\vee}_2\ra, &V_{-i}=\la f_1+i
f_2, f^{\vee}_1 + i f^{\vee}_2\ra,\\
&V_{\omega}=\la f_1, f^{\vee}_2\ra, &V_{\omega^2}=\la f_2,
f^{\vee}_1\ra.
\end{align}
\end{subequations}
 The proof is divided into two cases.

\noindent{\bf Case 1:} $p\equiv 1 \mod 3$. In the decomposition
(\ref{e2:supp}b), the operator $A$ commutes with $F$, it maps
$V_{\omega}$ isomorphically onto $V_{\omega^2}$, and $V_{\omega}$
and $V_{\omega^2}$ are $F$-invariant by Proposition 7.1. This
implies that
\[
char(V,F) = char(V_{\omega}, F)\cdot char(V_{\omega^2}, F) =
char(V_{\omega}, F)^2. \] Therefore $char(V_{\omega}, F)$ is
uniquely determined by $char(V, F)$. Let $\mu$ and $\nu$ be the
two roots of $char(V_{\omega}, F)$; then $\mu \ne \nu$ by
assumption.

As $F$ commutes with $A$, the operator $F$ also fixes $V_{\pm i}$.
We may assume that $\mu$ is a root of $char(V_i, F)$ and call the
other root $\lambda$. Let $v_1=a(f_1-i f_2) +
b(f^{\vee}_1-if^{\vee}_2)$ be an eigenvector of $F$ on $V_i$ with
eigenvalue $\mu$. Since $\mu$ is also an eigenvalue of $F$ on
$V_\omega$ and $V_{\omega^2}$, using the explicit bases of these
two spaces,  we may express $v_1$ as a sum of two
$\mu$-eigenvectors $\a f_1 + \b f^{\vee}_2\in V_{\omega}$ and $\c
f_2 + \d f^{\vee}_1\in V_{\omega^2}$ :
\[
a(f_1-i f_2) + b(f^{\vee}_1-if^{\vee}_2)=(\a f_1 + \b f^{\vee}_2)
+ (\c f_2 + \d f^{\vee}_1).
\]
Without loss of generality, we may assume $a = 1$. By comparing
the coefficients of both sides, we conclude that
$f_1-ibf^{\vee}_2$ and  $f_2+ibf^{\vee}_1$ form a basis of the
$\mu$-eigenspace of $F$ and $ v_1 =(f_1-ibf^{\vee}_2)
-i(f_2+ibf^{\vee}_1)$. Note that $v_2 = (f_1-ibf^{\vee}_2) +
i(f_2+ibf^{\vee}_1) = (f_1 + i f_2) - b (f^{\vee}_1+if^{\vee}_2)$
is a $\mu$-eigenvector of $F$ on $V_{-i}$. This shows that the
other eigenvalue $\lambda$ on $V_i$ is equal to $\nu$. Therefore
\[
char(V_i, F) = char(V_{-i}, F) = char(V_{\omega}, F) =
char(V_{\omega^2}, F).
\]

On the 2-adic side, we see this from the space $\widetilde W$.
More precisely, the space $\widetilde W_{\pm i}$ is the space
attached to the newform $g_{\pm}$, and $g_-$ is the twist of $g_+$
by the quadratic character $\chi_{-3}$ as noted in \S 6. For $p
\equiv 1 \mod 3$, we have $\chi_{-3}(p) = 1$ and hence
$char(\widetilde W_i, Frob_p) = char(\widetilde W_{-i}, Frob_p)$,
and consequently
\[
char(\widetilde W_i, Frob_p) = char(W_i, F_p) = char(V_i, F),
\] as desired.

\noindent{\bf Case 2:} $p\equiv 2 \mod 3$. In this case, by
Proposition 7.1, neither $F$ nor $F_p$ commutes with the operator
$\xi$, but $FA=AF$ and $F_pA=AF_p$ do. Therefore we use $FA$ and
$F_pA$ instead. The same analysis as in Case 1 yields
\[
char(V_{i}, FA) = char(V_{-i}, FA) = char(V_{\omega}, FA) =
char(V_{\omega^2}, FA).
\]
On the 2-adic side, as explained in Case 1, the action of $Frob_p$
on $W_{-i}$ is twisted by
 $\chi_{-3}(p) = -1$ of the action of $Frob_p$ on $W_i$, and hence
\[
char(\widetilde W_i, Frob_pH_{27}) = char(\widetilde W_{-i},
Frob_pH_{27}),
\] or equivalently,
\[
char(W_i, F_pA) = char(W_{-i}, F_pA).
\]
Applying (4.4.1) of \cite{Sch85} to the matrix $A$, we obtain
\[
char(V, FA) = char(W, F_pA),
\] hence
\[
char(V_i, FA)^2 = char(W_i, F_pA)^2,
\]
from which it follows that
\[
char(V_i, FA) = char(W_i, F_pA).
\] Since the action of $A$ on both spaces is multiplication by $i$, this yields
the desired equality
\[
char(V_i, F) = char(W_i, F_p).
\]
\end{proof}

The above theorem allows us to conclude that $f_{+}$ and $g_{+}$
(resp. $f_{-}$ and $g_{-}$) satisfy the Atkin-Swinnerton-Dyer
congruence relations at $p \ne 2, 3$ by the same argument as
Scholl's proof of Theorem 1.1 presented in section 5 of
\cite{Sch85}. This completes the proof of Theorem 1.2.
\medskip

\begin{rem} The same congruence relations at $p$ can be concluded from
Remark 5.8 of \cite{Sch85} provided that the action of $F$ is
ordinary.
\end{rem}

\bigskip

\section{The Atkin-Swinnerton-Dyer congruence relations and elliptic modular K3
surfaces over $\Q$} In this section we derive similar results for
the Atkin-Swinnerton-Dyer congruence relations arising from K3
surfaces over $\Q$.

First consider an explicit example.  Let $\G_2$ denote the group
associated to the algebraic equation (\ref{eqn_deg2}). A similar
discussion shows that this is a noncongruence subgroup. The space
$S_3(\G_2)$ of weight 3 cusp forms for $\G_2$ is 1-dimensional and
it is generated by
$$h_2=\sqrt{E_1E_2}=q^{1/10}-\frac{3^2}{2}q^{3/10}+\frac{3^3}{2^3}q^{5/10}+\frac{3\cdot
7^2}{2^4}q^{7/10}
  -\frac{3^2\cdot 7\cdot 19}{2^7}q^{9/10}+O(q^{10/10})$$ where $E_1$ and $E_2$
are given by (\ref{eqn_E1}) and (\ref{eqn_E2}) respectively. It is
clear that the  Fourier expansion of $h_2$ at the cusp $\infty$
has
 coefficients in the ring $\Z[1/2]$. Let
$g_2$ be a level 16 newform with the  first few terms of its
Fourier expansion (provided by Stein's data base \cite{stein}) as
$$g_2=q - 6q^5 + 9q^9 + 10q^{13} - 30q^{17} + 11q^{25} + 42q^{29} + O(q^{32})
$$
 Using Serre's method, we can show that the Atkin-Swinnerton-Dyer congruence
relation holds for $h_2$ and $g_2$. In fact, a more general result
regarding this situation can be proved.

 Recall that a {\it K3 surface} $S$ is a simply connected compact complex
surface with
 trivial canonical bundle. Its Hodge diamond, as the one defined
 in section \ref{some_elliptic_surface}, is
 \[
\begin{split}
&1\\
0\quad &\ \quad 0\\
1\quad\quad &20\quad\quad 1\\
0\quad &\ \quad 0\\
&1
\end{split}
\]

Hence as a compact complex surface, its Picard number $\rho(S)\le
h^{1,1}=20$. The cohomology group $\Lambda=H^2(S,\Z)$ is a rank 22
free $\Z$-module, called a K3 lattice. As a lattice, $\Lambda$ is
unimodular due to the Poincar\'e duality, even by Wu's formula,
with signature (3, 19) by the Hodge index theorem. Its
N\'eron-Severi group $NS(S)$, defined as in section
\ref{some_elliptic_surface}, is a sublattice of $\Lambda$ of
signature (1, $\rho(S)-1$) again by the Hodge index theorem. An
{\it elliptic surface} $\pi: S \rightarrow C$ is a two dimensional
complex variety over the base curve $C$ such that every fiber
$\pi^{-1}(t)$ is a smooth genus one curve except for finitely many
points $t$ in $C$. A compact elliptic surface is called an {\it
elliptic modular surface} if its monodromy group $\G_S$ is a
finite index subgroup of $SL_2(\Z)$ and $-I\notin \G_S$
~\cite{Shi72}.
\bigskip

\begin{theorem}
 Let $S$ be an elliptic modular K3 surface defined over $\Q$ with
 $\G_S$ being the associated modular group. Let $f_{\G_S}$ be a
nonzero $M$-integral form in the 1-dimensional space $S_3(\G_S)$
for some integer $M$.  Then there is a weight 3 cusp form
$g_{\G_S}$ with integral Fourier coefficients for some congruence
subgroup such that
 $f_{\G_S}$and $g_{\G_S}$ satisfy the Atkin-Swinnerton-Dyer congruence
 relations.
\end{theorem}
\begin{proof}
Since $S$ is an elliptic modular K3 surfaces, by Shioda's result
~\cite{Shi72}, its Picard number $\rho(S)=20$. Further, the work
 of Shioda and Inose \cite{s-i} on K3 surfaces with Picard number 20
shows that the Hasse-Weil $L$-function attached to $S$ contains a
factor $L(s, \chi^2)$, where $\chi$ is a Grossencharacter of some
imaginary quadratic
 extension of $\Q$ associated to an elliptic curve over $\Q$
with complex multiplications, arising from the K3 lattice of $S$.
 Combining the analytic behavior of
$L(s, \chi^2)$ proved by Hecke and the converse theorem for $GL_2$
proved by Weil, we know that $L(s, \chi^2)$ is also the
$L$-function attached to a weight 3 cuspidal newform $h_{\G_S}$
with integral coefficients for a congruence subgroup.

Let $\rho_l(S)$ be the $l$-adic representation of
$\rm{Gal}(\overline\Q/\Q)$ associated to $h_{\G_S}$, which exists
for almost all primes $l$. In particular, this representation is
isomorphic to $\rho_l^*$ defined in section 5 induced by the
action of $\rm{Gal}(\overline\Q/\Q)$ on $W_l$. This isomorphism of
representations can be seen explicitly via the characteristic
polynomials of Frobenius elements.  Due to a trace formula of
Monsky \cite{monsky71}, for almost all prime $p$, the Euler
$p$-factor $P_{21,p}(p^{-s})$ of the $L$-function of $h_{\G_S}$
appears as one part of the characteristic polynomial of the
Frobenius endomorphism $F_p$ acting on the crystalline cohomology
$H^2_{cris}(\CS_p/\Z_p)$, where $\CS_p$ is the N\'eron minimal
model of $S$ over $\F_p$. In particular, $F_{21,p}(p^{-s})$
corresponds to the characteristic polynomial of $F_p$ acting on
the orthogonal complement of the N\'eron-Severi group $NS(\CS_p
\otimes {\bar \F_p})\otimes \Q$ in $H^2_{cris}(\CS_p/\Z_p)\otimes
\Q$, where ${\bar \F_p}$ is the algebraic closure of $\F_p$  (for
details see \cite[\S 12]{b-s}). By the Lefschetz fixed point
theorem, the coefficients of $F_{21,p}(p^{-s})$ can be calculated
by  the same trace formula \cite[\S 3]{Sch88} for the Frobenius
endomorphism on $Frob_p$ on the $l$-adic cohomology $W_l$.

Denote by $\rho_l(\G_S)$ the $l$-adic representation of the Galois
group of $\Q$ associated to the space $S_3(\G_S)$  constructed by
Scholl in \cite{Sch85}. By Proposition 5.1, $\rho_l(\G_S)$ and
$\rho_l(S)$ are isomorphic up to a quadratic twist. Let $\phi$ be
a character of $\rm{Gal}(\overline\Q/\Q)$ of order at most two
such that $\rho_l(\G_S)$ is isomorphic to $\rho_l(S) \otimes
\phi$. By class field theory, we may regard $\phi$ as a Dirichlet
character of order at most two. Let $g_{\G_S}$ be the newform
having the same eigenvalues as the twist of $h_{\G_S}$ by $\phi$.
Then $g_{\G_S}$ also has integral Fourier coefficients and
$\rho_l(\G_S)$ is isomorphic to the $l$-adic representation
attached to $g_{\G_S}$. Applying Theorem 1.1 to $\rho_l(\G_S)$ and
$f_{\G_S}$, we conclude that $f_{\G_S}$ and $g_{\G_S}$ satisfy the
Atkin-Swinnerton-Dyer relation. \end{proof}

\end{document}